\documentclass{scrartcl}
\usepackage{graphicx,amsmath,amssymb}
\usepackage{algorithm,algorithmic}
\usepackage{appendix}
\usepackage{xcolor}
\usepackage{listings}
\usepackage[numbers]{natbib}
\usepackage{url}
\newtheorem{example}{Example}
\newcommand{\inprod}[2]{\left\langle #1, \, #2 \right\rangle}

\title{A Study on the Algorithm and Implementation of SDPT3}
\author{
    Naoki Ito
    \thanks{napinoco@gmail.com}
}

\date{\today}

\begin{document}

\maketitle
\begin{abstract}
This technical report presents a comprehensive study of SDPT3 \cite{toh1999,tutuncu2003,Toh2012}, a widely used open-source MATLAB solver for semidefinite-quadratic-linear programming, which is based on the interior-point method.
It includes a self-contained and consistent description of the algorithm, with mathematical notation carefully aligned with the implementation.
The aim is to offer a clear and structured reference for researchers and developers seeking to understand or build upon the implementation of SDPT3.
\end{abstract}

\tableofcontents

\section{Notations}\label{sec:notation}

For a positive integer $n$, let $\mathbb{R}^n$ denote the $n$-dimensional real vector space.
In this paper, we treat elements $x$ of $\mathbb{R}^n$ as column vectors, with the $i$-th component denoted by $x_i$.

For positive integers $m$ and $n$, let $\mathbb{R}^{m \times n}$ denote the space of $m \times n$ real matrices.
The $(i,j)$ component of $x \in \mathbb{R}^{m \times n}$ is denoted by $x_{ij}$.
For both vectors and matrices, the transpose of $x$ is denoted by $x^T$.
If $x \in \mathbb{R}^{n \times n}$ is invertible, its inverse is denoted by $x^{-1}$, and the transpose of the inverse by $x^{-T}$.
The trace of a matrix $x$ is denoted by $\operatorname{trace}(x)$, and the determinant by $\operatorname{det}(x)$.
For a vector $x \in \mathbb{R}^n$, we denote by $\operatorname{diag}(x) \in \mathbb{R}^{n \times n}$ the diagonal matrix with diagonal entries given by the components of $x$:
\[
\operatorname{diag}(x) = \begin{pmatrix}
x_1 & 0 & \cdots & 0 \\
0 & x_2 & \cdots & 0 \\
\vdots & \vdots & \ddots & \vdots \\
0 & 0 & \cdots & x_n
\end{pmatrix}.
\]
The identity matrix of size $n \times n$ is denoted by $I_n$ or simply $I$ when the dimension is clear from context.
Additionally, we define the diagonal matrix $J$ as follows:
\[
    J = \begin{pmatrix}
        1 & 0_{n-1}^T \\
        0_{n-1} & -I_{n-1}
    \end{pmatrix} \in \mathbb{R}^{n \times n},
\]
where $0_{n-1} \in \mathbb{R}^{n-1}$.

Regarding the notation $\|a\|$, if $a$ is a matrix, it represents the Frobenius norm, and if $a$ is a vector, it represents the L2 norm (Euclidean norm).
Specifically,
\[
\|a\| = 
\begin{cases}
    \sqrt{\sum_{i=1}^m \sum_{j=1}^n a_{ij}^2} & \text{if } a \in \mathbb{R}^{m \times n},\\
    \sqrt{\sum_{i=1}^n a_i^2} & \text{if } a \in \mathbb{R}^n.
\end{cases}
\]

We define several sets as follows:
\begin{itemize}
    \item 
    \textbf{$n$-dimensional non-negative real cone:} \\
    \[
      \mathbb{R}^n_+ 
      = \{\,x \in \mathbb{R}^n \mid x \geq 0 \}.
    \]

    \item 
    \textbf{$n$-dimensional second-order cone:} \\
    \[
      \mathbb{Q}^n 
      = \left\{ \begin{pmatrix} x_0 \\ \bar{x} \end{pmatrix} \in \mathbb{R}^n 
         \;\middle|\; x_0 \in \mathbb{R}^1_+, \; \bar{x} \in \mathbb{R}^{n-1}, \; \|\bar{x}\| \le x_0 \right\}.
    \]
    Here, we define a non-negative real-valued function $\gamma : \mathbb{Q}^n \to \mathbb{R}^1_+$ as
    \[
      \gamma(x) = \sqrt{x^T J\, x}.
    \]

    \item 
    \textbf{Space of $n$-dimensional real symmetric matrices:} \\
    \[
      \mathbb{S}^n = \{\, x \in \mathbb{R}^{n \times n} \mid x = x^T \}.
    \]

    \item 
    \textbf{$n$-dimensional real positive semi-definite cone:} \\
    \[
      \mathbb{S}^n_+ 
      = \{\, x \in \mathbb{S}^n \mid a^T x\, a \ge 0 \;\; \forall a \in \mathbb{R}^n \}.
    \]
\end{itemize}

For $x \in \mathbb{R}^n$ and $\epsilon > 0$, the open ball of radius $\epsilon$ centered at $x$ is defined as
\[
  B(x, \epsilon) = \{ y \in \mathbb{R}^n \mid \|y - x\| < \epsilon \}.
\]
For a set $S \subseteq \mathbb{R}^n$, the affine hull of $S$ is defined as
\[
  \operatorname{aff}(S) = \left\{ \sum_{i=1}^k \lambda_i x_i \mid x_i \in S, \lambda_i \in \mathbb{R}, \sum_{i=1}^k \lambda_i = 1, k \in \mathbb{N} \right\}.
\]
The interior of a set $S$ is defined as
\[
  \operatorname{int}(S) 
  = \{\, x \in S \mid \exists\, \epsilon > 0 \quad \text{s.t.} \quad B(x, \epsilon) \subseteq S \},
\]
whereas the relative interior is defined as the interior within the affine hull of $S$:
\[
\operatorname{relint}(S)
=\{\,
x \in S
\mid \exists\, \epsilon > 0 \quad \text{s.t.} \quad \bigl(\operatorname{aff}(S)\cap B(x,\epsilon)\bigr)\,\subseteq\,S
\}.
\]
In this paper, we use $\operatorname{relint}$ instead of $\operatorname{int}$ to properly handle the singleton set $\{0\}$:
we have $\operatorname{relint}(\{0\}) = \{0\}$ while $\operatorname{int}(\{0\}) = \emptyset$.
For all proper cones considered in this paper ($\mathbb{R}^n_+$, $\mathbb{Q}^n$, $\mathbb{S}^n_+$), the two notions coincide.

\section{Problem Definition}

In this paper, we consider the following primal problem (P) and its corresponding dual problem (D):
\begin{equation*}
    \everymath{\displaystyle}
    \renewcommand{\arraystretch}{2.0}
    \text{(P)}~
    \left|
    \begin{array}{cl}
         \min_{x} & \displaystyle 
             \sum_{p \in P} \inprod{c^p}{x^p}_p 
             \;+\;\sum_{p \in P} \phi^p\bigl(x^p;\,\nu^p\bigr) \\[3pt]
         \text{s.t.} 
         & \displaystyle 
             \sum_{p \in P} \inprod{a^p_{k}}{x^p}_p 
             \;=\; b_k \quad (k= 1,2,\ldots,m), \\[3pt]
         & x^p \;\in\; \mathbb{K}^p \quad (p\in P),
    \end{array}
    \right.
    \qquad
    \text{(D)}~
    \left|
    \begin{array}{cl}
         \max_{y,z} & \displaystyle 
            \sum_{k=1}^m b_k\, y_k 
            \;+\;\sum_{p \in P} \phi^{p*}\bigl(z^p;\,\nu^p\bigr) \\[3pt]
         \text{s.t.} 
         & \displaystyle 
            c^p \;-\; \sum_{k=1}^m a^p_k\,y_k \;=\; z^p \quad (p\in P), \\[3pt]
         & y \;\in\; \mathbb{R}^m, \quad z^p \;\in\; (\mathbb{K}^p)^* \quad (p\in P)
    \end{array}
    \right.
\end{equation*}
where the following parameters are assumed to be given:
\begin{itemize}
    \item A positive integer $p_{\max}$ representing the number of cone blocks, and the index set $P = \{1, 2, \ldots, p_{\max}\}$.
    \item For each block $p \in P$, a positive integer $n_p$ representing the dimension of the variables.
    \item A positive integer $m$ representing the number of constraints.
    \item The cone $\mathbb{K}^p$ for each $p \in P$, which is one of $\mathbb{S}^{n_p}_+$, $\mathbb{Q}^{n_p}$, $\mathbb{R}^{n_p}_+$, or $\mathbb{R}^{n_p}$.
    \item Coefficients $c^p$ and $a^p_k$ for $k = 1, \ldots, m$ and $p \in P$, where:
      \begin{itemize}
          \item $c^p, a^p_k \in \mathbb{S}^{n_p}$ if $\mathbb{K}^p = \mathbb{S}^{n_p}_+$
          \item $c^p, a^p_k \in \mathbb{R}^{n_p}$ otherwise
      \end{itemize}
    \item A coefficient vector $b \in \mathbb{R}^m$.
    \item Non-negative parameters $\nu^p \ge 0$ for all $p \in P$.
\end{itemize}
The notation $(\mathbb{K}^p)^*$ denotes the dual cone of $\mathbb{K}^p$, 
i.e., $(\mathbb{K}^p)^* = \{ z^p \mid \inprod{x^p}{z^p}_p \geq 0 \text{ for all } x^p \in \mathbb{K}^p \}$, 
which in our case takes the following forms:
\[
(\mathbb{K}^p)^* = 
\begin{cases}
    \mathbb{S}_{+}^{n_p} & \text{if } \mathbb{K}^{p} = \mathbb{S}_{+}^{n_p},\\
    \mathbb{Q}^{n_p}     & \text{if } \mathbb{K}^{p} = \mathbb{Q}^{n_p},\\
    \mathbb{R}_{+}^{n_p} & \text{if } \mathbb{K}^{p} = \mathbb{R}_{+}^{n_p},\\
    \{0\}                & \text{if } \mathbb{K}^{p} = \mathbb{R}^{n_p}.
\end{cases}
\]
The inner product $\left\langle a, x \right\rangle_p$ is defined as:
\[
\left\langle a, x \right\rangle_p = 
\begin{cases}
    \operatorname{trace}(a^T x) 
    = \sum_{i=1}^{n_p} \sum_{j=1}^{n_p} a_{ij} x_{ij}, 
    & \text{if } \mathbb{K}^p = \mathbb{S}_{+}^{n_p},\\[6pt]
    a^T x = \sum_{i=1}^{n_p} a_i x_i, 
    & \text{otherwise}.
\end{cases}
\]
The function $\phi^p : \mathbb{K}^p \to \mathbb{R}_+$ is a barrier function defined as:
\[
\phi^p(x^p;\, \nu^p) =
\begin{cases}
    -\nu^p \log \det(x^p), & \text{if } \mathbb{K}^{p} = \mathbb{S}_{+}^{n_p},\\
    -\nu^p \log \gamma(x^p), & \text{if } \mathbb{K}^{p} = \mathbb{Q}^{n_p},\\
    -\sum_{i=1}^{n_p} \nu^p \log x^p_i, & \text{if } \mathbb{K}^{p} = \mathbb{R}_{+}^{n_p},\\
    0, & \text{if } \mathbb{K}^{p} = \mathbb{R}^{n_p},
\end{cases}
\]
where $\gamma(x^p) = \sqrt{(x^p)^T J x^p}$ as introduced in Section~\ref{sec:notation}.
The notation $\phi^{p*}$ denotes the convex conjugate of $\phi^p$, 
i.e., $\phi^{p*}(z^p; \nu^p) = \sup_{x^p \in \mathbb{K}^p} \{\langle z^p, x^p \rangle_p - \phi^p(x^p; \nu^p)\}$, 
which in our case takes the following forms:
\[
\phi^{p*}(z^p;\, \nu^p) =
\begin{cases}
    \nu^p \log \det(z^p) + n_p \nu^p (1 - \log \nu^p), 
    & \text{if } \mathbb{K}^{p} = \mathbb{S}_{+}^{n_p},\\[4pt]
    \nu^p \log \gamma(z^p) + \nu^p (1 - \log \nu^p), 
    & \text{if } \mathbb{K}^{p} = \mathbb{Q}^{n_p},\\[4pt]
    \sum_{i=1}^{n_p} \left( \nu^p \log z^p_i + \nu^p (1 - \log \nu^p) \right), 
    & \text{if } \mathbb{K}^{p} = \mathbb{R}_{+}^{n_p},\\[3pt]
    0, & \text{if } \mathbb{K}^{p} = \mathbb{R}^{n_p}.
\end{cases}
\]

When $\nu^p = 0$ for all $p \in P$, problems (P) and (D) reduce to the standard form of conic linear programming.
By considering the case $\nu^p > 0$, we can handle a broader class of problems, including log-determinant optimization problems.
\section{Infeasible Primal-Dual Path-Following Interior-Point Method}
\label{sec:infeasible_IPM}
This section presents an infeasible primal-dual path-following interior-point method that converges to a feasible solution from an arbitrary infeasible starting point.
Such infeasible-start methods, which are employed in SDPT3, are particularly important in practice for two reasons: finding an initial feasible point can be computationally expensive even when feasibility is known, and numerical errors in problem data may make theoretically feasible points practically infeasible.

For simplicity of notation, we define two linear mappings 
$\mathcal{A}^p : \mathbb{K}^p \to \mathbb{R}^m$ and 
$(\mathcal{A}^p)^T : \mathbb{R}^m \to \mathbb{K}^p$ 
as follows:
\[
  \mathcal{A}^p x^p
  := 
  \begin{pmatrix}
      \left\langle a^p_1, x^p \right\rangle_p \\
      \left\langle a^p_2, x^p \right\rangle_p \\
      \vdots \\
      \left\langle a^p_m, x^p \right\rangle_p
  \end{pmatrix},
  \qquad
  (\mathcal{A}^p)^T y
  :=
  \sum_{k=1}^m a^p_k y_k.
\]
Then, the primal problem (P) and the dual problem (D) can be written as follows:
\[
  \everymath{\displaystyle}
  \renewcommand{\arraystretch}{2.0}
  \text{(P)}~
  \left|
  \begin{array}{cl}
      \min_{x} & \displaystyle
       \sum_{p \in P} \left\langle c^p, x^p \right\rangle_p - \sum_{p \in P} \phi^p(x^p;\, \nu^p) \\[3pt]
      \text{s.t.} & \displaystyle
       \sum_{p \in P} \mathcal{A}^p x^p = b,\\[3pt]
      & x^p \in \mathbb{K}^p \quad (p \in P),
  \end{array}
  \right.
  \qquad
  \text{(D)}~
  \left|
  \begin{array}{cl}
      \max_{y,z} & \displaystyle
       \sum_{k=1}^m b_k y_k 
          + \sum_{p \in P} \phi^{p*}(z^p;\, \nu^p) \\[3pt]
      \text{s.t.} & \displaystyle
       c^p - (\mathcal{A}^p)^T y = z^p \quad (p \in P),\\[3pt]
      & y \in \mathbb{R}^m, \quad z^p \in (\mathbb{K}^p)^* \quad (p \in P).
  \end{array}
  \right.
\]

We define $P^{\text{u}} := \{\, p \in P \mid \mathbb{K}^p = \mathbb{R}^{n_p} \,\}$.\footnote{%
  The superscript $\text{u}$ stands for "unrestricted," meaning $\mathbb{K}^p = \mathbb{R}^{n_p}$.
}
The KKT conditions for the optimality of problems (P) and (D) are:
\begin{equation}
    \everymath{\displaystyle}
    \renewcommand{\arraystretch}{2.5}
    \left\{
    \begin{array}{ll}
        \sum_{p \in P} \mathcal{A}^p x^p - b = 0, & \\[-4pt]
        (\mathcal{A}^p)^T y + z^p - c^p = 0, & \quad (p \in P),\\[-4pt]
        x^p \circ z^p - \nu^p e^p = 0, & \quad (p \in P \setminus P^{\text{u}}),\\[-4pt]
        x^p \in \mathbb{K}^p,\; y \in \mathbb{R}^m,\; z^p \in (\mathbb{K}^p)^*, & \quad (p \in P).
    \end{array}
    \right.
    \label{eq:KKTcond}
\end{equation}
In the third equation, the bilinear mapping $x^p \circ z^p$ is defined for each block $p \in P \setminus P^{\text{u}}$ as
\[
  x^p \circ z^p = 
  \begin{cases}
    \frac{1}{2} \left( x^p (z^p)^T + z^p (x^p)^T \right), 
      & \text{if } \mathbb{K}^p = \mathbb{S}^{n_p}_+,\\[4pt]
    \left( (x^p)^T z^p;\; x^p_0 \bar{z}^p + z^p_0 \bar{x}^p \right),
      & \text{if } \mathbb{K}^p = \mathbb{Q}^{n_p},\\[4pt]
    \operatorname{diag}(x^p) z^p,
      & \text{if } \mathbb{K}^p = \mathbb{R}^{n_p}_+,
  \end{cases}
\]
where $e^p$ is the identity element for this operator:
\[
  e^p = 
  \begin{cases}
    I, & \mathbb{K}^p = \mathbb{S}^{n_p}_+,\\[3pt]
    (1,\,0,\ldots,0)^T, & \mathbb{K}^p = \mathbb{Q}^{n_p},\\[3pt]
    (1,\,1,\ldots,1)^T, & \mathbb{K}^p = \mathbb{R}^{n_p}_+ 
  \end{cases}
\]
The bilinear operator $\circ$ is called the Jordan product on $\mathbb{K}^p$.
\textbf{
Note that in the case of $\mathbb{S}^{n_p}_+$, it is usually defined as $x^p \circ z^p = \frac{1}{2}(x^p z^p + z^p x^p)$, 
but in this paper, we extend this operator to non-symmetric matrices and use the notation $x^p \circ z^p = \frac{1}{2} \left( x^p (z^p)^T + z^p (x^p)^T \right)$ for the sake of notational simplicity.
While this extension of the operator $\circ$ to non-symmetric matrices does not satisfy the definition of Jordan algebra \cite{Faraut1994}, it allows for unified and simplified notation (especially in equations like \eqref{eq:NewtonKKT}).
}
Now, consider replacing each $\nu^p$ ($p \in P$) in the complementarity conditions with a common positive parameter $\mu > 0$. 
Then, under mild regularity conditions, the perturbed KKT system \eqref{eq:KKTcond} provides a unique solution $\big( x(\mu),\, y(\mu),\, z(\mu) \big)$ within the feasible region. 
This solution forms a differentiable trajectory with respect to $\mu$, called the central path $T$:
\[
  T = \left\{ \big( x(\mu),\, y(\mu),\, z(\mu) \big) \mid \mu > 0 \right\}.
\]
In practice, larger values of $\mu$ tend to yield better numerical conditioning. Therefore, the path-following method starts with a large value of $\mu$ and gradually decreases it toward the target values $\nu^p$, tracking the central path to find the optimal solution of the original problem.

In the following sections, we will detail the main components of this algorithm, including the definition and computation of the search direction, the choice of step size, and convergence criteria.

\subsection{Search Directions} \label{sec:direction}
\subsubsection{Framework of the Search Direction}
We assume that a current iterate $(x,\,y,\,z) \in \operatorname{int}(\mathbb{K}) \times \mathbb{R}^m \times \operatorname{relint}(\mathbb{K}^*)$ is given, where $\mathbb{K} := \mathbb{K}^1 \times \cdots \times \mathbb{K}^{p_{\max}}$ and $\mathbb{K}^* := (\mathbb{K}^1)^* \times \cdots \times (\mathbb{K}^{p_{\max}})^*$, with $x = (x^1, \ldots, x^{p_{\max}})$ and $z = (z^1, \ldots, z^{p_{\max}})$.
\footnote{
  We use $\operatorname{relint}$ to handle the case where $\mathbb{K}^p = \mathbb{R}^{n_p}$, for which the dual cone $(\mathbb{K}^p)^* = \{0\}$ satisfies $\operatorname{relint}(\{0\}) = \{0\} \neq \operatorname{int}(\{0\}) = \emptyset$. 
  For $\mathbb{K}^p = \mathbb{S}^{n_p}_+$, $\mathbb{Q}^{n_p}$, or $\mathbb{R}^{n_p}_+$, we have $\operatorname{relint}((\mathbb{K}^p)^*) = \operatorname{int}((\mathbb{K}^p)^*)$.
}

To compute the next iterate along the central path, we apply Newton's method to the perturbed KKT system. 
Introducing the parameter $\sigma \in [0,1)$ and the scaling matrices $G^p$ for $p \in P$, the search direction $(\Delta x,\, \Delta y,\, \Delta z)$ is obtained as the solution to the following equations:
\begin{equation}
    \renewcommand{\arraystretch}{2.5}
    \left\{
    \begin{array}{rll}
         \sum_{p \in P} \mathcal{A}^p \Delta x^p & = R_{\text{prim}} := b - \sum_{p \in P} \mathcal{A}^p x^p &  \\
         (\mathcal{A}^p)^T \Delta y + \Delta z^p & = R_{\text{dual}}^p := c^p - z^p - (\mathcal{A}^p)^T y & (p \in P) \\
         \mathcal{E}^p \Delta x^p + \mathcal{F}^p \Delta z^p & = R_{\text{comp}}^p := \max\{\sigma \mu, \nu^p\} e^p - \big(G^p x^p\big) \circ \big((G^p)^{-1} z^p\big) & (p \in P \setminus P^{\text{u}})
    \end{array}
    \right.
    \label{eq:NewtonKKT}
\end{equation}
where 
\begin{equation}
  \mu := \frac{\sum_{p \in P^0} \left\langle x^p, z^p \right\rangle}
              {\sum_{p \in P^0} n^p}
  \label{eq:mu}
\end{equation}
with $P^0 := \{p\in P \mid \nu^p=0\}$, and the linear mappings $\mathcal{E}^p : \mathbb{K}^p \to \mathbb{K}^p$ and $\mathcal{F}^p : \mathbb{K}^p \to \mathbb{K}^p$ are defined as
\[
  \mathcal{E}^p \Delta x^p
    := \big(G^p \Delta x^p\big) \circ \big((G^p)^{-1} z^p\big),
  \quad
  \mathcal{F}^p \Delta z^p
    := \big(G^p x^p\big) \circ \big((G^p)^{-1} \Delta z^p\big).
\]

Equation \eqref{eq:NewtonKKT} represents a linearized system derived from Newton's method applied to the perturbed KKT conditions \eqref{eq:KKTcond}. 
The first two equations directly correspond to the primal and dual feasibility conditions, while the third equation linearizes the complementarity conditions with the parameter $\sigma$ and scaling matrices $G^p$.
The parameter $\sigma \in [0,1)$ controls the target point on the central path, with $\sigma \mu$ serving as the target duality gap.
The scaling matrices $G^p$ provide variable transformations that affect the convergence properties. 
Different choices of $G^p$ lead to distinct search directions with varying computational costs and numerical stability, as detailed in the following subsections.

\subsubsection{AHO Direction}
Consider first the simplest scaling choice where $G^p = I$. 
This yields the Alizadeh-Haeberly-Overton (AHO) search direction \cite{Alizadeh1998}, which is widely used in semidefinite programming and second-order cone programming. 
The AHO search direction exhibits good numerical stability but incurs higher computational cost per iteration compared to scaled variants. 
For linear and unrestricted blocks ($\mathbb{K}^p = \mathbb{R}^{n_p}_+$ and $\mathbb{K}^p = \mathbb{R}^{n_p}$), the identity scaling $G^p = I$ remains the standard choice.

\subsubsection{HKM Direction}
For semidefinite and second-order cone blocks ($\mathbb{K}^p = \mathbb{S}^{n_p}_+$ or $\mathbb{K}^p = \mathbb{Q}^{n_p}$), consider the following scaling matrices $G^p$.

\paragraph{Semidefinite cone: \(\mathbb{K}^p = \mathbb{S}^{n_p}_+\).}
\[
  G^p := (z^p)^{\tfrac12}.
\]
In this case, $G^p e^p (G^p)^T = z^p$ holds.

\paragraph{Second-order cone: \(\mathbb{K}^p = \mathbb{Q}^{n_p}\).}
\[
  \omega^p := \gamma(z^p), 
  \quad
  t^p := \frac{1}{\gamma(z^p)}\, z^p,
  \quad
  G^p :=
  \omega^p
  \begin{pmatrix}
    t^p_0 & (\bar{t}^p)^T \\
    \bar{t}^p & I + \frac{1}{1+t^p_0}\,\bar{t}^p(\bar{t}^p)^T
  \end{pmatrix}.
\]
In this case, $G^p e^p = z^p$ holds.

\medskip 

The search direction obtained by solving equation \eqref{eq:NewtonKKT} with these scaling matrices $G^p$ is called the HKM search direction.\footnote{
More precisely, this approach was proposed by Helmberg--Rendl--Vanderbei--Wolkowicz~\cite{Helmberg1996} / Kojima--Shindoh--Hara~\cite{Kojima1997} / Monteiro~\cite{Monteiro1997}, and is also referred to as the HRVW/KSH/M direction.
} 
It is known for its numerical stability and effectiveness in exploiting problem sparsity, making it particularly suitable for large-scale semidefinite programs.

\subsubsection{NT Direction}
For semidefinite and second-order cone blocks ($\mathbb{K}^p = \mathbb{S}^{n_p}_+$ or $\mathbb{K}^p = \mathbb{Q}^{n_p}$), consider the following scaling matrices $G^p$.

\paragraph{Semidefinite cone: \(\mathbb{K}^p = \mathbb{S}^{n_p}_+\).}
\[
  G^p 
  := 
    \Bigl( (x^p)^{\frac12} \bigl( (x^p)^{\frac12} z^p (x^p)^{\frac12} \bigr)^{-\frac12} (x^p)^{\frac12} \Bigr)^{-\frac12}.
\]
In this case, $(G^p)^{-1} z^p (G^p)^{-T} = G^p x^p (G^p)^T$ holds. 

\paragraph{Second-order cone: \(\mathbb{K}^p = \mathbb{Q}^{n_p}\).}
\begin{equation}
    \omega^p := \sqrt{\frac{\gamma(z^p)}{\gamma(x^p)}}, 
    \quad 
    \xi^p 
    := \begin{pmatrix} \xi^p_0 \\ \bar{\xi}^p \end{pmatrix} 
    = \begin{pmatrix}
        \frac{1}{\omega^p} z^p_0 + \omega^p x^p_0 \\
        \frac{1}{\omega^p} \bar{z}^p - \omega^p \bar{x}^p \\
    \end{pmatrix},
    \quad
    t^p := \frac{1}{\gamma(\xi^p)}\xi^p
    \label{eq:scaling_mat_NT_socp_aux}
\end{equation}
\begin{equation}
    G^p := \omega^p \begin{pmatrix}
        t^p_0 & (\bar{t}^p)^T \\
        \bar{t}^p & I+\frac{1}{1 + t^p_0} \bar{t}^p(\bar{t}^p)^T
    \end{pmatrix}
    \label{eq:scaling_mat_NT_socp}
\end{equation}
In this case, $(G^p)^{-1} z^p = G^p x^p$ holds.

\medskip 

The search direction obtained by solving equation \eqref{eq:NewtonKKT} with these scaling matrices $G^p$ is called the Nesterov--Todd (NT) search direction \cite{Nesterov1997,todd1998}. 
The NT search direction is based on the theory of self-concordant barriers and is known for its strong theoretical guarantees of global convergence.

\medskip 

In SDPT3, \texttt{HKM} and \texttt{NT} options are available for semidefinite and second-order cone blocks ($\mathbb{K}^p = \mathbb{S}^{n_p}_+$ and $\mathbb{K}^p = \mathbb{Q}^{n_p}$), while the identity scaling $G^p = I$ is always used for linear and unrestricted blocks ($\mathbb{K}^p = \mathbb{R}^{n_p}_+$ and $\mathbb{K}^p = \mathbb{R}^{n_p}$), regardless of the selected search direction.

\subsection{Reduction of the equation}
In this section, we reduce the search direction equation \eqref{eq:NewtonKKT} into a more manageable form. 
Specifically, we eliminate $\Delta x^p$ and $\Delta z^p$ to obtain an equation in terms of $\Delta y$, resulting in the so-called Schur complement system.

First, by solving $(\mathcal{A}^p)^T \Delta y + \Delta z^p = R_{\text{dual}}^p$ and $\mathcal{E}^p \Delta x^p + \mathcal{F}^p \Delta z^p = R_{\text{comp}}^p$ of \eqref{eq:NewtonKKT} for $\Delta z^p$ and $\Delta x^p$ respectively, we obtain the equation in terms of $\Delta y$:
\begin{equation}
    \everymath{\displaystyle}
    \renewcommand{\arraystretch}{1.5}
    \left\{
    \begin{array}{rll}
    \Delta z^p &= R_{\text{dual}}^p - \mathcal{A}^p(\Delta y)  & (p\in P)\\
    \Delta x^p &= (\mathcal{E}^p)^{-1}R_{\text{comp}}^p - \mathcal{H}^p(\Delta z^p) \\
               &= (\mathcal{E}^p)^{-1}R_{\text{comp}}^p - \mathcal{H}^p(R_{\text{dual}}^p - \mathcal{A}^p(\Delta y))  & (p\in P\setminus P^{\text{u}})
   \end{array}
   \right.
   \label{eq:sol_x_z}
\end{equation}
where $\mathcal{H}^p := (\mathcal{E}^p)^{-1}\mathcal{F}^p$. 
The existence of $(\mathcal{E}^p)^{-1}$ is non-trivial, but it is known to exist if $G^p$ is positive definite \cite{todd1998}. 
Since all $G^p$ introduced in Section~\ref{sec:direction} are positive definite, $(\mathcal{E}^p)^{-1}$ exists.
Substituting \eqref{eq:sol_x_z} into the first equation of \eqref{eq:NewtonKKT}, we obtain:
\begin{equation}
    \left\{
    \begin{aligned}
        \sum_{p \in P\setminus P^{\text{u}}} \mathcal{A}^p\mathcal{H}^p(\mathcal{A}^p)^T\Delta y + \sum_{p \in P^{\text{u}}} \mathcal{A}^p(\Delta x^p) 
            &= R_{\text{prim}} - \sum_{p \in P\setminus P^{\text{u}}} \mathcal{A}^p \big( (\mathcal{E}^p)^{-1}R_{\text{comp}}^p - \mathcal{H}^p R_{\text{dual}}^p \big) \\
        (\mathcal{A}^p)^T \Delta y 
            &= R^p_{\text{dual}} \qquad (p\in P^{\text{u}})
    \end{aligned}
    \right.
    \label{eq:Schur_complement}
\end{equation}

\medskip

Next, we consider reducing \eqref{eq:Schur_complement} to matrix representations. 
Let the matrix $M^p \in \mathbb{S}^{m}$ satisfy
\[
    \mathcal{A}^p\mathcal{H}^p(\mathcal{A}^p)^T \Delta y= M^p \Delta y,
\]
and define:
\begin{align*}
    M &= \sum_{p \in P \setminus P^{\text{u}}} M^p \\
    h &= R_{\text{prim}} - \sum_{p \in P \setminus P^{\text{u}}} \mathcal{A}^p\big((\mathcal{E}^p)^{-1}R_{\text{comp}}^p - \mathcal{H}^p R_{\text{dual}}^p\big)\\
    A^p &= \begin{pmatrix}
        (a^p_1)^T\\
        (a^p_2)^T\\
        \vdots\\
        (a^p_m)^T
    \end{pmatrix} \in \mathbb{R}^{m\times n_p} \quad (p\in P)\\
    A^{\text{u}} &= [ A^p ~ (p\in P^{\text{u}})\text{ concatenated horizontally} ]\\
    R^{\text{u}}_{\text{dual}} &= [ R^p_{\text{dual}} ~ (p\in P^{\text{u}})\text{ concatenated vertically} ]\\
    \Delta x^{\text{u}} &= [ \Delta x^p ~ (p\in P^{\text{u}})\text{ concatenated vertically} ]
\end{align*}
Then, the equation \eqref{eq:Schur_complement} can be represented in matrix form as:
\begin{equation}
    \underbrace{\left(\begin{array}{cc}
        M   & A^{\text{u}} \\
        A^{\text{u}} & O
    \end{array}\right)}_{\mathcal{M}}
    \left(\begin{array}{c}
        \Delta y   \\
        \Delta x^{\text{u}} 
    \end{array}\right) 
    = 
    \left(\begin{array}{c}
         h  \\
         R_{\text{dual}}^{\text{u}} 
    \end{array}
    \right)
    \label{eq:Schur_complement_Mat}
\end{equation}

\medskip

In the following, we describe the specific computation of the matrix $M$ and vector $h$ in the Schur complement system \eqref{eq:Schur_complement_Mat}, namely the calculation methods for $(\mathcal{E}^p)^{-1}R^p_{\text{comp}}$, $\mathcal{H}^p R^p_{\text{dual}}$, and the matrix $M^p$ for both HKM and NT search directions. 
The detailed derivations are omitted here; a concise guide can be found in Appendix~\ref{sec:guide_for_dir_eq}. 

For both HKM and NT search directions, we have:
\begin{equation}
    (\mathcal{E}^p)^{-1}R^p_{\text{comp}} = \max\{\sigma\mu, \nu^p\}(z^p)^{-J} - x^p \qquad (p \in P)
    \label{eq:Einv_Rcomp}
\end{equation}
where $(z^p)^{-J}$ is the inverse element of $z^p$ with respect to the Jordan product $\circ$:
\[
    (z^p)^{-J} = \begin{cases}
        (z^p)^{-1} & \text{if } \mathbb{K}^p = \mathbb{S}^{n_p}_+, \\
        \frac{1}{(\gamma(z^p))^2} J z^p & \text{if } \mathbb{K}^p = \mathbb{Q}^{n_p}, \\
        (1/z_1, \ldots, 1/z_{n_p})^T & \text{if } \mathbb{K}^p = \mathbb{R}^{n_p}_+.
    \end{cases}
\]
\medskip
When the search direction option is set to \texttt{HKM}, SDPT3 uses:
\begin{align}
    \mathcal{H}^p R^p_{\text{dual}} &= \begin{cases}
        \frac{1}{2}(x^p R^p_{\text{dual}} (z^p)^{-1} + (z^p)^{-1} R^p_{\text{dual}} x^p) & \text{if } \mathbb{K}^p = \mathbb{S}^{n_p}_+,\\
        - \Big( (x^p)^T J (z^p)^{-J} \Big) J R^p_{\text{dual}} + \inprod{(z^p)^{-J}}{R^p_{\text{dual}}} x^p + \inprod{R^p_{\text{dual}}}{x^p} (z^p)^{-J} & \text{if } \mathbb{K}^p = \mathbb{Q}^{n_p},\\
        \operatorname{diag}(x^p) \operatorname{diag}(z^p)^{-1} \Delta z^p & \text{if } \mathbb{K}^p = \mathbb{R}^{n_p}_+
    \end{cases} \label{eq:HKM_HRd}\\
    M^p &= \begin{cases} 
        \text{a matrix whose $(k,\ell)$-elements are given by } \inprod{a^p_k}{x^p a^p_\ell (z^p)^{-1}} & \text{if } \mathbb{K}^p = \mathbb{S}^{n_p}_+, \\
        -\Big( (x^p)^T J (z^p)^{-J} \Big) A^p J (A^p)^T + (A^p x^p)(A^p (z^p)^{-J})^T + (A^p (z^p)^{-J})(A^p x^p)^T& \text{if } \mathbb{K}^p = \mathbb{Q}^{n_p}, \\
        A^p \operatorname{diag}(x^p) \operatorname{diag}(z^p)^{-1} (A^p)^T & \text{if } \mathbb{K}^p = \mathbb{R}^{n_p}_+.
    \end{cases}
\end{align}
\medskip
When the search direction option is set to \texttt{NT}, SDPT3 uses:
\begin{align}
    \mathcal{H}^p R^p_{\text{dual}} &= \begin{cases}
        W^p R^p_{\text{dual}} W^p & \text{if } \mathbb{K}^p = \mathbb{S}^{n_p}_+, \\
        \frac{1}{(\omega^p)^2} \Big(-J R^p_{\text{dual}} + 2\inprod{R^p_{\text{dual}}}{t^p}t^p\Big) & \text{if } \mathbb{K}^p = \mathbb{Q}^{n_p},\\
        \operatorname{diag}(x^p) \operatorname{diag}(z^p)^{-1} R^p_{\text{dual}} & \text{if } \mathbb{K}^p = \mathbb{R}^{n_p}_+
    \end{cases}\\
    M^p &= \begin{cases}
        \text{a matrix whose $(k,\ell)$-elements are given by } \inprod{a^p_k}{W^p\, a^p_\ell\, W^p} & \text{if } \mathbb{K}^p = \mathbb{S}^{n_p}_+,\\
        \frac{1}{(\omega^p)^2}\Big(-A^p J (A^p)^T + 2 (A^p t^p)(A^p t^p)^T \Big) & \text{if } \mathbb{K}^p = \mathbb{Q}^{n_p}, \\
        A^p \operatorname{diag}(x^p) \operatorname{diag}(z^p)^{-1} (A^p)^T & \text{if } \mathbb{K}^p = \mathbb{R}^{n_p}_+.
    \end{cases} \label{eq:NT_M}
\end{align}
where $W^p = (G^p)^{-2} = (x^p)^{-\frac{1}{2}}((x^p)^{\frac{1}{2}} z^p (x^p)^{\frac{1}{2}})^{\frac{1}{2}} (x^p)^{-\frac{1}{2}}$, and $\omega^p$, $t^p$ are defined as in \eqref{eq:scaling_mat_NT_socp_aux}.

Although $W^p$ appears to be very complex, it can be calculated using the following procedure \cite{todd1998}:
\begin{enumerate}
    \item First, perform the Cholesky decomposition of $z^p$ to obtain the upper triangular matrix $U$, i.e., $z^p=U^TU$.
    \item Next, perform the eigenvalue decomposition of $U x^p U$ to obtain the orthogonal matrix $V$ and the diagonal matrix $\Lambda$ with eigenvalues on the diagonal, i.e., $U x^p U = V \Lambda V^T$.
    \item Now, let $S=\Lambda^\frac{1}{4}(U^{-1}V)^T$, then $W^p=S^T S$.
\end{enumerate}

\subsection{Solving the reduced equation} \label{sec:solve_reduced_eq}
The matrix $\mathcal{M}$ in \eqref{eq:Schur_complement_Mat} is often very ill-conditioned, making the direct solution of \eqref{eq:Schur_complement_Mat} numerically unstable. 
To address this issue, SDPT3 employs iterative methods with preconditioning. 
Specifically, SDPT3 \cite{toh1999} uses either the Symmetric Quasi-Minimal Residual (SQMR) method \cite{Freund1994} or the Biconjugate Gradient Stabilized (BiCGSTAB) method \cite{vanderVorst1992,Golub2013}.

In both iterative methods, it is necessary to repeatedly calculate the product of the preconditioner (approximation of $\mathcal{M}^{-1}$) and a vector.
Two approaches can be considered for efficient and accurate computation of this matrix-vector product:
The first approach is to perform LU decomposition of the entire matrix $\mathcal{M}$ in advance.
This allows the product of $\mathcal{M}^{-1}$ and a vector to be calculated by forward and backward substitution, which is numerically more stable than explicitly computing the inverse of the ill-conditioned matrix $\mathcal{M}$.

If $M$ is positive definite, another approach can be taken. 
Using the Schur complement $S := (A^{\text{u}})^T M^{-1} A^{\text{u}} - O$, we have:
\[
\mathcal{M}^{-1}=\begin{pmatrix}
    M^{-1} + M^{-1} A^{\text{u}} S^{-1} (A^{\text{u}})^T M^{-1} & -M^{-1} A^{\text{u}} S^{-1} \\
    -S^{-1} (A^{\text{u}})^T M^{-1} & S^{-1}
\end{pmatrix}
\]
Thus, we obtain:
\[
    \mathcal{M}^{-1}\begin{pmatrix}u \\ v \end{pmatrix} = \begin{pmatrix} \hat{u} - M^{-1} A^{\text{u}} \hat{v} \\ \hat{v} \end{pmatrix}
\]
where $\hat{u} = M^{-1} u$ and $\hat{v} = S^{-1}\big((A^{\text{u}})^T \hat{u} - v \big)$.
If the Cholesky decomposition of $M$ and the LU decomposition of $S$ are calculated in advance, the calculation of the product of $\mathcal{M}^{-1}$ and a vector reduces to three forward-backward substitutions for $M$ and two forward-backward substitutions for $S$.
This approach is computationally more efficient for two reasons: $M$ can use Cholesky decomposition (which is faster than LU decomposition), and $S$ is typically much smaller than $\mathcal{M}$.
However, while this method is efficient from a computational standpoint, it has the disadvantage that $S$ tends to become extremely ill-conditioned.

In the implementation of SDPT3, the second method is tried first, and if the Cholesky decomposition of $M$ fails or if $S=LU$ is found to be extremely ill-conditioned (specifically, if the ratio of the maximum to minimum diagonal elements of $U$ is greater than $10^{30}$), the first method of LU decomposition of the entire $\mathcal{M}$ is adopted.

The choice of iterative method in SDPT3 depends on the preconditioning approach used.
SDPT3 primarily uses the numerically stable SQMR method, but SQMR requires the preconditioner to be symmetric.
When Cholesky-based preconditioning is feasible, SQMR can be safely used.
However, when the matrix is ill-conditioned and LU decomposition must be used, the LU factorization may not accurately preserve the symmetry of the original matrix. 
Therefore, SDPT3 switches to BiCGSTAB for safety, as it does not require the preconditioner to be symmetric.\footnote{
  However, there is an implementation inconsistency: in the predictor step (\texttt{linsysolve.m}, called by \texttt{HKMpred.m} and \texttt{NTpred.m}), SDPT3 version 4.0 still uses the SQMR method even with LU preconditioning. 
  The empirical success of SDPT3 across diverse problem instances suggests that SQMR remains robust to the small numerical asymmetries that may arise in LU-based preconditioning.
}

Both approaches can be further optimized by exploiting the sparsity of the problem. 
Details are discussed in Section~\ref{sec:exploit_sparsity_socp_lp}.

\subsection{Step size} \label{sec:step_size}
In this section, we explain the method for calculating the step sizes $\alpha_P$, $\alpha_D \in [0, 1]$ used to determine the next iteration point 
\[
  (x^+,\,y^+,\,z^+) 
  \;=\; 
  \bigl(x + \alpha_P\,\Delta x,\;\; y + \alpha_D\,\Delta y,\;\; z + \alpha_D\,\Delta z\bigr).
\]
The step sizes must be chosen such that the new iterates satisfy
\[
  x^+   
  \;\in\; \operatorname{int}\bigl(\mathbb{K}\bigr),
  \quad
  z^+   
  \;\in\; \operatorname{relint}\bigl((\mathbb{K})^*\bigr).
\]
To ensure this, we first compute the maximum feasible step lengths
\[
  \alpha_x  := \sup \bigl\{\,\alpha \in [0,1] \mid x + \alpha \,\Delta x \in \mathbb{K}\bigr\},
  \quad
  \alpha_z  := \sup \bigl\{\,\alpha \in [0,1] \mid z + \alpha \,\Delta z \in (\mathbb{K})^*\bigr\}
\]
and then use a constant $\gamma\in(0,1)$ (e.g., $\gamma=0.99$) to obtain the actual step sizes
\[
  \alpha_P 
    = \gamma \,\alpha_x, 
  \quad
  \alpha_D 
    = \gamma \,\alpha_z.
\]

\subsubsection{Computation of $\alpha_x$}
Let \[
  \alpha^p_x:=\sup\{\alpha \geq 0 \mid x^p + \alpha \Delta x^p \in \mathbb{K}^p\} \quad (p\in P)
\]
then $\alpha_x = \min\{1, \min\{\alpha^p_x \mid p \in P\}\}$ holds.
Below we describe the computation of the maximum feasible step length $\alpha^p_x$ for each cone type.

\paragraph{Semidefinite cone: $\mathbb{K}^p = \mathbb{S}^{n_p}_+$.}
Using the Cholesky decomposition $x^p=LL^T$ of $x^p$,
\[\alpha^p_x = \sup\{\alpha \geq 0 \mid I + \alpha L^{-1} \Delta x^p L^{-T} \in \mathbb{K}^p\}\]
holds, thus if $\lambda_{\max}$ is the largest eigenvalue of $-L^{-1} \Delta x^p L^{-T}$,
\[\alpha^p_x = \begin{cases}
    1/\lambda_{\max} & \text{if } \lambda_{\max} > 0, \\
    +\infty & \text{otherwise}
\end{cases}\]
is obtained.
$\lambda_{\max}$ can be efficiently calculated with good accuracy using methods such as the Lanczos method \cite{Golub2013}.

\paragraph{Second-order cone: $\mathbb{K}^p = \mathbb{Q}^{n_p}$.}
Let
\begin{align*}
  f^p(\alpha)
  :&= ( x^p + \alpha \Delta x^p )^T J ( x^p + \alpha \Delta x^p )\\
   &=  ( x^p_0 + \alpha \Delta x^p_0 )^2 - \bigl\|( \bar{x}^p + \alpha \Delta \bar{x}^p )\bigr\|^2.
\end{align*}
Then we have
\begin{align*}
    \alpha_x^p 
    &= \sup\{\alpha \geq 0 \mid x^p + \alpha \Delta x^p \in \mathbb{Q}^{n_p}\}\\
    &= \sup\{\alpha \geq 0 \mid f^p(\alpha) \geq 0, ~x^p_0 + \alpha \Delta x^p_0 \geq 0\}.
\end{align*}

Consider the non-negative roots of the quadratic equation $f^p(\alpha)=0$. 
Expanding $f^p(\alpha)$, we get 
\begin{align*}
    f^p(\alpha) = \alpha^2\underbrace{(\Delta x^p)^T J (\Delta x^p)}_a + 2 \alpha \underbrace{(x^p)^T J (\Delta x^p)}_b + \underbrace{(x^p)^T J x^p}_c,
\end{align*}
where the discriminant is $d:=b^2-ac$.
Since $x^p\in \operatorname{int}(\mathbb{Q}^{n_p})$ is assumed, we have $c = (x^p_0)^2 - \|\bar{x}^p\|^2 > 0$.
Hence the quadratic equation $f^p(\alpha) = 0$ has positive solutions only in the following three cases:
\begin{enumerate}
    \item When $a<0$, $\alpha=\frac{-b-\sqrt{d}}{a}$ is the unique positive solution.
    \item When $a=0$ and $b<0$, $\alpha=-\frac{c}{2b}$ is the unique positive solution.
    \item When $a>0$, $b<0$, and $d\geq 0$, $\alpha=\frac{-b-\sqrt{d}}{a}$ is the smallest positive solution.
\end{enumerate}
Furthermore, these solutions always satisfy $x^p_0+\alpha \Delta x^p_0\geq 0$.
(\textbf{Proof:} If $\Delta x^p_0 \geq 0$, it is trivially satisfied. 
Consider the case $\Delta x^p_0 < 0$. 
Note that $f^p(0)=c>0$ and $f^p(-\frac{x^p_0}{\Delta x^p_0}) = - \|(x^p -\frac{x^p_0}{\Delta x^p_0} \Delta x^p)\|^2 \leq 0$.
Therefore, the smallest positive solution of $f^p(\alpha) = 0$ exists in the interval $\big(0, \frac{x^p_0}{-\Delta x^p_0}\big]$, 
and throughout this interval $x^p_0+\alpha \Delta x^p_0\geq 0$ holds.)

Therefore,
\begin{equation*}
    \alpha^p_x = \begin{cases}
       \frac{-b - \sqrt{d}}{a} & \text{if } (a < 0) \text{ or } (a > 0 \text{ and } b < 0 \text{ and } d \geq 0),\\
       -\frac{c}{2b} & \text{if } a = 0 \text{ and } b < 0,\\
       +\infty & \text{otherwise}.
    \end{cases}
\end{equation*}

\paragraph{Linear cone: $\mathbb{K}^p = \mathbb{R}^{n_p}_+$.}
For each $i = 1, 2, \ldots, n_p$, define
\[
t_i = \begin{cases}
    -x^p_i / \Delta x^p_i & \text{if } \Delta x^p_i < 0, \\
    +\infty & \text{otherwise}.
\end{cases}
\]
Then
\[
    \alpha^p_x = \min\{t_i \mid i = 1, 2, \ldots, n_p\}.
\]

\paragraph{Unrestricted variables: $\mathbb{K}^p = \mathbb{R}^{n_p}$.}
In this case, $\alpha^p_x = +\infty$ since there are no constraints on $x^p$.

\subsubsection{Computation of $\alpha_z$}
To compute $\alpha_z$, we first calculate $\alpha^p_z = \sup\{\alpha \geq 0 \mid z^p + \alpha\Delta z^p \in (\mathbb{K}^p)^*\}$ for each block $p \in P$, then take $\alpha_z=\min\{1, \min\{\alpha^p_z \mid p \in P\}\}$.
Since $\mathbb{K}^p=\mathbb{S}^{n_p}_+,\mathbb{Q}^{n_p},\mathbb{R}^{n_p}_+$ are self-dual cones ($\mathbb{K}^p = (\mathbb{K}^p)^*$), $\alpha^p_z$ can be calculated in the same way as $\alpha^p_x$.
For $\mathbb{K}^p=\mathbb{R}^{n_p}$, since $(\mathbb{K}^p)^*=\{0\}^{n_p}$, we have $\Delta z^p = 0$ and thus $\alpha^p_z=+\infty$.

\subsection{Initial Points}
\label{sec:initial_points}
While the algorithms described in this paper can start from infeasible initial points, the choice of initial points significantly affects both convergence speed and numerical stability.
It is reported that numerical calculations tend to become unstable when given initial points with extremely small or large norms \cite{toh1999}.
Therefore, it is desirable to provide initial points with scales comparable to the expected solutions of problems (P) and (D).
SDPT3 employs the following empirically effective initialization: $y = 0$, and for each $p\in P$,
\[
    x^p = \begin{cases}
        \zeta^p\, e^p, & \text{if } p \in P \setminus P^{\text{u}},\\
        0, & \text{if } p\in P^{\text{u}},
    \end{cases}
    \quad
    \quad 
    z^p = \begin{cases}
        \eta^p\, e^p, & \text{if } p \in P \setminus P^{\text{u}},\\
        0, & \text{if } p\in P^{\text{u}},
    \end{cases}
\]
where
\begin{align*}
    \zeta^p 
    &= \max\Bigl\{
       10,\;\sqrt{n^p},\;\;\theta^p \max_{1 \le k \le m}\bigl\{\frac{1 + |b_k|}{1 + \|a^p_k\|}\bigr\}
      \Bigr\},
    \quad
    \quad
    \theta^p = \begin{cases}
        n^p, & \mathbb{K}^p = \mathbb{S}^{n_p}_+,\\
        \sqrt{n^p}, & \mathbb{K}^p = \mathbb{Q}^{n_p},\\
        1, & \text{otherwise},
    \end{cases}
    \\[6pt]
    \eta^p 
    &= \max\Bigl\{
       10,\;\sqrt{n^p},\;\max\{\|a^p_1\|,\ldots,\|a^p_m\|,\;\|c^p\|\}
      \Bigr\}.
\end{align*}

\subsection{Stopping Criteria}
\label{sec:stopping_criteria}
In this section, we present the stopping criteria for terminating the iterations of SDPT3.
The algorithm terminates iterations when the predefined number of iterations or accuracy goals are met, or when infeasibility or numerical difficulties become apparent.

Specifically, we first define the following quantities to measure the dual gap and infeasibility:
\[
  \mathtt{gap}
  := \sum_{p\in P}
       \Bigl(\inprod{x^p}{z^p}_p + \bigl(\phi^p(x^p) - \phi^{p*}(z^p)\bigr)\Bigr),
\]
\[
  \mathtt{relgap}
  := \frac{\mathtt{gap}}
           {\,1 \;+\;\Bigl|\sum_{p\in P}\inprod{c^p}{x^p}_p\Bigr|
                 \;+\;|\,b^T y|\,},
\]
\[
  \mathtt{pinfeas}
  := \frac{\|\,R_{\text{prim}}\|}{\,1 + \|b\|\,},
  \quad
  \mathtt{dinfeas}
  := \frac{\sum_{p\in P}\|\,R^p_{\text{dual}}\|}
           {\,1 + \sum_{p\in P}\|\,c^p\|\,}.
\]
where $R_{\text{prim}}$ and $R_{\text{dual}}^p$ are defined in equation \eqref{eq:NewtonKKT}.
These represent the primal and dual infeasibility residuals, respectively.
Note that if $\mathtt{pinfeas} = 0$ and $\mathtt{dinfeas} = 0$,
\[
  \sum_{p\in P}\inprod{x^p}{z^p}_p 
  \;=\;
  \sum_{p\in P}\inprod{x^p}{c^p - (\mathcal{A}^p)^T y}_p
  \;=\; \sum_{p\in P}\inprod{x^p}{c^p}_p \;-\; b^T y,
\]
holds, indicating that $\mathtt{gap}$ can be used as an indicator of the dual gap.

Under these definitions, the iterations are terminated when any of the following conditions are met:
\begin{enumerate}
    \item The number of iterations reaches the upper limit.
    \item 
      $\displaystyle
      \max\{\,\mathtt{relgap},\;\mathtt{pinfeas},\;\mathtt{dinfeas}\}
      < \varepsilon
      $, i.e., the desired accuracy is achieved.
    \item 
      $\displaystyle
        \frac{\,b^T y\,}
              {\,\sum_{p\in P}\|\,(\mathcal{A}^p)^T y + z^p\|\!}
      > \kappa
      $
      indicating that the primal problem (P) is likely infeasible.
    \item 
      $\displaystyle
      -\,\frac{\inprod{c}{x}}
              {\bigl\|\sum_{p\in P}\mathcal{A}^p x^p\bigr\|}
      > \kappa
      $
      indicating that the dual problem (D) is likely infeasible.
    \item Numerical errors occur:
      \begin{itemize}
          \item Failure in Cholesky decomposition of $x^p$ or $z^p$.
          \item Preconditioned iterative solvers (SQMR or BiCGSTAB) fail to converge.
          \item $\mathtt{gap}$ diverges, indicating erratic behavior.
      \end{itemize}
\end{enumerate}
In practice, it is also important to terminate iterations when $\mathtt{relgap}$, $\mathtt{pinfeas}$, and $\mathtt{dinfeas}$ take relatively small values and there is little improvement in the last few iterations.
In SDPT3 \cite{toh1999}, various heuristic termination conditions are implemented in \texttt{sqlpcheckconvg.m}.

\section{Predictor-Corrector Method}
In this section, we introduce the predictor-corrector method.\footnote{The name predictor-corrector originates from numerical methods for ordinary differential equations.}
This method is widely adopted in many software packages as a practically efficient approach for solving optimization problems.
The predictor-corrector method was originally proposed by Mehrotra \cite{Mehrotra1992}.
Wright \cite{Wright1997} later provided a clear and theoretically elegant presentation of the method, which has become a standard reference.
Here, we introduce the formulation by Toh et al. \cite{toh1999}, which is equivalent to the formulation presented by Wright but differs slightly in its computational procedure.

The predictor-corrector method computes the search direction $(\Delta x, \Delta y, \Delta z)$ in two stages.
In the first stage (predictor step), we compute an affine scaling direction $(\delta x, \delta y, \delta z)$ by setting $\sigma = 0$, which aims to reduce the duality gap aggressively.
Based on the progress achieved by this predictor direction, we determine an appropriate centering parameter $\sigma$, and in the second stage (corrector step), we compute the final search direction $(\Delta x, \Delta y, \Delta z)$ that balances progress toward optimality with proximity to the central path.

\paragraph{Predictor Step}
First, we set $\sigma=0$ to compute the affine scaling direction.
In other words, we replace $R^p_{\text{comp}}$ in \eqref{eq:NewtonKKT} with
\[R^{p,\text{pred}}_{\text{comp}}=\nu^p \, e^p - (G^p x^p) \circ ((G^p)^{-1} z^p).\]
Solving \eqref{eq:NewtonKKT} with this modification yields the predictor direction $(\delta x, \delta y, \delta z)$.
Next, we compute trial step sizes $\alpha_P, \alpha_D$ for this predictor direction using the method described in Section~\ref{sec:step_size}, without actually updating the variables.
These step sizes are used to estimate the potential reduction in the duality gap and to determine an appropriate centering parameter $\sigma$.
Specifically, the centering parameter is computed as
\[
   \sigma=\min\left\{1, \left(\frac{\inprod{x + \alpha_P \delta x}{z + \alpha_D \delta z}}{\inprod{x}{z}}\right)^\psi\right\},
\]
where $\psi \ge 1$ is an exponent parameter.
Empirical studies suggest that $\psi \in \{2,3,4\}$ works well.
SDPT3 adaptively selects $\psi$ based on a default parameter $\hat{\psi}=3$:
\[
\psi = \begin{cases}
    \max\{\hat{\psi}, 3 \min(\alpha_P, \alpha_D)^2\} & \text{if} ~ \mu > 10^{-6}, \\
    \max\{1, \min\{\hat{\psi}, 3 \min(\alpha_P, \alpha_D)^2\}\} & \text{otherwise}.
\end{cases}
\]
Recall that $\mu$ is defined as in \eqref{eq:mu}.

\paragraph{Corrector Step}
With the computed centering parameter $\sigma$, we now solve for the final search direction.
The complementarity residual in \eqref{eq:NewtonKKT} is replaced with
\[R^{p,\text{corr}}_{\text{comp}}=\max\{\sigma \mu, \nu^p\} e^p - \big(G^p x^p\big) \circ \big((G^p)^{-1} z^p\big) - \big(G^p \delta x^p\big)\circ\big((G^p)^{-1} \delta z^p\big),\]
where the last term is the second-order correction from the predictor step.
Solving \eqref{eq:NewtonKKT} with this modification yields the final search direction $(\Delta x, \Delta y, \Delta z)$.
We then compute the step sizes $\beta_P, \beta_D$ for this search direction using the method in Section~\ref{sec:step_size}, and update the solution as
\[(x^+, y^+, z^+) = (x, y, z) + (\beta_P \Delta x, \beta_D \Delta y, \beta_D \Delta z).\]

\bigskip
A key advantage of the predictor-corrector method is its computational efficiency.
Between the predictor and corrector steps, only the complementarity residual $R^p_{\text{comp}}$ changes, which affects only the right-hand side vector $h$ in the Schur complement system \eqref{eq:Schur_complement_Mat}.
The coefficient matrix $\mathcal{M}$ and its factorization (LU or Cholesky) computed in the predictor step can be reused in the corrector step.
Furthermore, the right-hand side vector can be efficiently updated.
Letting $h_{\text{pred}}$ and $h_{\text{corr}}$ denote the vectors for the predictor and corrector steps respectively, we have
\[h_{\text{corr}}=h_{\text{pred}} + \sum_{p\in P\setminus P^{\text{u}}} \mathcal{A}^p (\mathcal{E}^p)^{-1} \big((G^p \delta x^p) \circ ((G^p)^{-1} \delta z^p) \big).\]
Thus, although the predictor-corrector method requires solving the linear system twice, the additional computational cost is minimal due to these reuse strategies.

\section{Sparsity Exploitation Technique} \label{sec:exploit_sparsity}
\subsection{Semidefinite cone} \label{sec:exploit_sparsity_sdp}
The dominant computational cost in the algorithms described so far arises from the calculation of $M^p$ when $\mathbb{K}^p = \mathbb{S}^{n_p}_+$. 
Naively computing all $O(m^2)$ elements of $M^p$ is computationally expensive since each element requires matrix multiplications and inner products, making it essential to exploit the sparsity of $a^p_k$ for practical efficiency. 
SDPT3 implements a modified version of the technique proposed by Fujisawa et al.~\cite{Fujisawa1997} as described in \cite{tutuncu2001}. 
Here, we present this approach with further reorganization for clarity of exposition.

\medskip

Recall that, in the HKM direction, we have
\[
  M^p_{ij}
  = \inprod{a^p_i}{x^p\,a^p_j\,(z^p)^{-1}}
\]
Let $j$ be fixed. Observe that:
\begin{enumerate}
\item[(i)] Once we compute $G = x^p\,a^p_j\,(z^p)^{-1}$, it can be reused to calculate $M^p_{ij} = \langle a^p_i, G \rangle$ for all $i$.
\item[(ii)] Since $M^p$ is symmetric, we only need to compute the upper triangular part (i.e., $M^p_{ij}$ for $i \leq j$).
\item[(iii)] For computing $M^p_{ij}$, we only need the elements of $G$ corresponding to the non-zero positions of $a^p_i$ (where $i \leq j$), rather than all elements of $G$.
\end{enumerate}

To leverage property (iii) efficiently, we introduce a permutation $\sigma$ of indices $1,\ldots,m$ such that $f_{\sigma(k)}$ is in ascending order, where $f_k$ denotes the number of non-zero elements in $a^p_k$.
For a fixed $j$, define the index set
\[
  I
  := \bigl\{(\alpha,\beta)\mid (a^p_{\sigma(i)})_{\alpha\beta}\neq 0
         \text{ for some } i=1,\ldots,j \bigr\}
\]
representing all matrix positions needed for inner products with $a^p_{\sigma(i)}$ ($i \leq j$).
Consequently, when computing $G = x^p\,a^p_{\sigma(j)}\,(z^p)^{-1}$, we need only determine the element $G_{\alpha\beta}$ where $(\alpha,\beta) \in I$. 
By processing matrices in order of increasing density, the size of $I$ tends to grow gradually, enabling significant computational savings in the early iterations (i.e., for small $j$).

We now describe three methods (F1--F3) for computing these required elements, evaluating their computational costs by counting the number of multiplications:
\begin{enumerate}
\item[\textbf{F1:}] First, calculate $F=a^p_{\sigma(j)}\,(z^p)^{-1}$, which requires $n_p f_{\sigma(j)}$ multiplications.
           Then, calculate $G=x^p F$, which requires $n_p^3$ multiplications.
           Thus, a total of $n_p f_{\sigma(j)} + n_p^3$ multiplications are required.
\item[\textbf{F2:}] First, calculate $F=a^p_{\sigma(j)}\,(z^p)^{-1}$, which requires $n_p f_{\sigma(j)}$ multiplications. Then compute:
  \[
    G_{\alpha\beta}=
    \begin{cases}
     \sum_{\gamma=1}^{n_p} (x^p)_{\alpha\gamma} \, F_{\gamma\beta}, & \text{if } (\alpha,\beta)\in I,\\
     0, & \text{otherwise},
    \end{cases}
  \]
  This computation requires $n_p$ multiplications for each of the $|I|$ elements.
  Thus, a total of $n_p (f_{\sigma(j)} + |I|)$ multiplications are required.
\item[\textbf{F3:}] Compute each element of $G$ directly as follows:
  \[
    G_{\alpha\beta}=
    \begin{cases}
      \sum_{(\gamma,\delta): (a^p_{\sigma(j)})_{\gamma\delta} \neq 0}\,
        (x^p)_{\alpha\gamma} \, (a^p_{\sigma(j)})_{\gamma\delta} \, (z^p)^{-1}_{\delta\beta},
       & \text{if } (\alpha,\beta) \in I,\\
      0, & \text{otherwise}.
    \end{cases}
  \]
  Since each non-zero element of $a^p_{\sigma(j)}$ contributes to the computation, this requires $2f_{\sigma(j)}$ multiplications for each of the $|I|$ elements.
  Thus, a total of $2f_{\sigma(j)}|I|$ multiplications are required.
\end{enumerate}
F1 is straightforward but computationally expensive, though it benefits from highly optimized BLAS routines \cite{Dongarra1990_BLAS} for dense matrix multiplication.
F2 and F3 compute only the necessary elements $(\alpha,\beta) \in I$, which can lead to significant savings when $|I| \ll n_p^2$, despite some overhead from less cache-friendly memory access patterns compared to BLAS operations.
Specifically, F2 is efficient when the set $I$ is small relative to the matrix size (i.e., $|I| \ll n_p^2$), while F3 excels when the matrix $a^p_{\sigma(j)}$ is highly sparse (i.e., $f_{\sigma(j)} \ll n_p$).
In SDPT3, the choice between methods F1--F3 for each $j$ is determined during a preprocessing phase, based on the sparsity patterns and expected computational costs.

While we have presented these methods for the HKM direction, the same sparsity exploitation techniques apply to the NT direction.

\subsection{Second-order and linear cones} \label{sec:exploit_sparsity_socp_lp}
When solving equation \eqref{eq:Schur_complement_Mat} using the methods described in Section~\ref{sec:solve_reduced_eq},
if the matrices $\mathcal{M}$ and $M$ are sparse, we can speed up matrix-vector multiplications within the SQMR or BiCGSTAB methods and utilize highly optimized Cholesky decomposition and LU decomposition routines such as CHOLMOD \cite{Chen2008_CHOLMOD} and UMFPACK \cite{Davis2004_UMFPACK}.

Unfortunately, in practical large-scale problems, the matrices $\mathcal{M}$ and $M$ are frequently dense. 
However, $M$ can often be expressed as a sparse positive (semi)definite matrix $M_{\mathrm{sparse}}$ plus a low-rank perturbation.
This structure arises because most constraint vectors $a^p_k$ are sparse, with only a few being dense.
Below, we describe a sparsity exploitation technique that can be utilized in such cases.

Assume that for any $p\in P\setminus P^{\text{u}}$, the matrix $M^p$ can be decomposed as
\begin{equation}
  M^p = M^p_{\mathrm{sparse}} + U^p\, D^p\, (U^p)^T
  \label{eq:low_rank_perturbation}
\end{equation}
where $M^p_{\mathrm{sparse}} \in \mathbb{S}^m_+$ is a sparse positive semidefinite matrix, 
$U^p \in \mathbb{R}^{m\times n^p_+}$ with $n^p_+ \ll m$ (low rank), 
and $D^p\in \mathbb{S}^{n^p_+}$.
By aggregating these decompositions, we define:
\begin{align*}
  M_{\mathrm{sparse}} &:= \sum_{p\in P\setminus P^{\text{u}}} M^p_{\mathrm{sparse}}, \\
  U &:= [U^p ~ (p\in P) \text{ concatenated horizontally}], \\
  D &:= [D^p ~ (p\in P) \text{ concatenated block diagonally}].
\end{align*}
Using this decomposition, we can transform the original system \eqref{eq:Schur_complement_Mat} 
into an augmented sparse system. Specifically, the following linear equation
\begin{equation}
  \begin{pmatrix}
    M_{\mathrm{sparse}} & A^{\text{u}} & U \\
    (A^{\text{u}})^T & O & O \\
    U^T & O & -D^{-1}
  \end{pmatrix}
  \begin{pmatrix}
    \Delta y \\
    \Delta x^{\text{u}} \\
    \lambda
  \end{pmatrix}
  =
  \begin{pmatrix}
    h \\
    R^{\text{u}}_{\text{dual}} \\
    0
  \end{pmatrix}
  \label{eq:Schur_complement_Mat_aug}
\end{equation}
is equivalent to \eqref{eq:Schur_complement_Mat} with the auxiliary variable $\lambda = D\,U^T\,\Delta y$.

The augmented coefficient matrix has dimension $m+\sum_{p\in P\setminus P^{\text{u}}} n^p_+$, 
which is larger than the original $m \times m$ system. 
However, it offers the computational advantage of being sparse, 
enabling the use of efficient sparse linear algebra routines.
This system \eqref{eq:Schur_complement_Mat_aug} can be solved using the methods 
described in Section~\ref{sec:solve_reduced_eq} with 
$A^{\text{u}}$ replaced by $[\,A^{\text{u}}, U]$ and the zero blocks appropriately replaced by 
$\begin{pmatrix} O & O \\ O & -D^{-1} \end{pmatrix}$.

Unfortunately, when $\mathbb{K}^p = \mathbb{S}^{n_p}_+$ for some $p \in P$, 
$M^p$ is often a dense matrix without exploitable structure, 
making it difficult to express in the form of \eqref{eq:low_rank_perturbation}.
However, when $\mathbb{K}^p = \mathbb{Q}^{n_p}$ or $\mathbb{K}^p = \mathbb{R}^{n_p}_+$, 
we can construct the decomposition \eqref{eq:low_rank_perturbation} 
and utilize the augmented system \eqref{eq:Schur_complement_Mat_aug}.

We now describe specific methods for constructing $M^p_{\mathrm{sparse}}$, $U^p$, and $D^p$ 
for these cases.
First, partition $A^p$ into sparse and dense columns:
$A^p_{\mathrm{sparse}}$ contains the sparse columns and 
$A^p_{\mathrm{dense}}$ contains the dense columns.
In SDPT3, this partition is based on the sparsity ratio of each column;
a column is classified as sparse if its ratio of non-zero elements 
is below a certain threshold (0.4 by default).
If $A^p_{\mathrm{dense}}$ is empty, we simply set $M^p_{\mathrm{sparse}} = M^p$, 
with $U^p$ and $D^p$ as empty matrices.
Otherwise, we construct the decomposition as follows.

\paragraph{Second-order cone:}
Consider the case $\mathbb{K}^p=\mathbb{Q}^{n_p}$ with the HKM direction.
Using the identity $-J = I - 2\,e^p(e^p)^T$, 
the expression for $M^p$ from Section~\ref{sec:direction} can be rewritten as
\[
  M^p 
  = \bigl((x^p)^T J (z^p)^{-J}\bigr)\, A^p(A^p)^T
    \;+\; u^p (v^p)^T
    \;+\; v^p (u^p)^T
    \;-\; 2\,((x^p)^T J (z^p)^{-J})\, k^p (k^p)^T,
\]
where $u^p := A^p x^p$, $v^p := A^p (z^p)^{-J}$, and $k^p := A^p e^p$.
When $A^p_{\mathrm{dense}}$ is non-empty, the vectors $u^p$ and $v^p$ become dense, 
making $M^p$ a dense matrix.

To exploit sparsity, SDPT3 decomposes $M^p$ by separating the contributions 
from sparse and dense columns. Specifically, it defines:
\[
  M^p_{\mathrm{sparse}}
    := \bigl((x^p)^T J (z^p)^{-J}\bigr)\,
       A^p_{\mathrm{sparse}}\,(A^p_{\mathrm{sparse}})^T,
\]
\[
  U^p
    := \Bigl(
       \sqrt{(x^p)^T J (z^p)^{-J}}\;A^p_{\mathrm{dense}},
       \quad u^p,
       \quad \gamma(z^p)^2\,v^p,
       \quad -\sqrt{2\,(x^p)^T J (z^p)^{-J}}\,k^p
    \Bigr),
\]
\[
  D^p
    := \begin{pmatrix}
         I & O & O & O \\
         O & 0 & 1/\gamma(z^p)^2 & 0 \\
         O & 1/\gamma(z^p)^2 & 0 & 0 \\
         O & 0 & 0 & -1
       \end{pmatrix}.
\]
This decomposition satisfies $M^p = M^p_{\mathrm{sparse}} + U^p D^p (U^p)^T$, 
where the inclusion of the $k^p$ term ensures that $M^p_{\mathrm{sparse}}$ 
remains positive semidefinite. 
Note that in the implementation, $D^p$ is never constructed; 
only $-(D^p)^{-1}$ is directly computed:
\[
  -(D^p)^{-1}
  = \begin{pmatrix}
      -I & O & O & O \\
      O & 0 & -\gamma(z^p)^2 & 0 \\
      O & -\gamma(z^p)^2 & 0 & 0 \\
      O & 0 & 0 & 1
    \end{pmatrix}
\]

For the NT direction, the expression for $M^p$ takes a different form.
From Section~\ref{sec:direction}, we have
\[
  M^p
  = \frac{1}{(\omega^p)^2} \Big(
      A^p (A^p)^T + 2 u^p (u^p)^T - 2k^p (k^p)^T
  \Big),
\]
where $u^p := A^p t^p$ and $k^p := A^p e^p$.
When $A^p_{\mathrm{dense}}$ is non-empty, the vector $u^p$ becomes dense, making $M^p$ a dense matrix.

Similar to the HKM direction, SDPT3 decomposes $M^p$ by separating the contributions from sparse and dense columns:
\[
  M^p_{\mathrm{sparse}}
    := \frac{1}{(\omega^p)^2} \, A^p_{\mathrm{sparse}}\,(A^p_{\mathrm{sparse}})^T,
\]
\[
  U^p
    := \Bigl(
      \frac{1}{\omega^p}\,A^p_{\mathrm{dense}},
       \;\; \sqrt{2}\,u^p,
       \;\; -\sqrt{2}\,k^p
    \Bigr),
\]
\[
  D^p
    := \begin{pmatrix}
         I & O & O \\
         O & 1/(\omega^p)^2 & 0  \\
         O & 0 & -1 
       \end{pmatrix}.
\]
This decomposition satisfies $M^p = M^p_{\mathrm{sparse}} + U^p D^p (U^p)^T$,
where again the inclusion of the $k^p$ term ensures that $M^p_{\mathrm{sparse}}$ 
remains positive semidefinite.
As with the HKM direction, $D^p$ is never constructed in the implementation;
only $-(D^p)^{-1}$ is directly computed:
\[
  -(D^p)^{-1}
  = \begin{pmatrix}
      -I & O & O  \\
      O & -(\omega^p)^2 & 0  \\
      O & 0 & 1
    \end{pmatrix}.
\]

\paragraph{Linear cone:}
Consider the case $\mathbb{K}^p = \mathbb{R}^{n_p}_+$. 
Unlike the second-order cone case, SDPT3 uses the same formula ($G^p=I$) 
for both \texttt{HKM} and \texttt{NT} options.
Recall from Section~\ref{sec:direction} that
\[
  M^p = A^p \operatorname{diag}(x^p) \operatorname{diag}(z^p)^{-1} (A^p)^T.
\]
When $A^p_{\mathrm{dense}}$ is non-empty, this matrix becomes dense.

To exploit sparsity, SDPT3 decomposes $M^p$ as follows.
Partition $x^p$ and $z^p$ according to the sparse/dense partition of $A^p$:
$x^p_{\mathrm{sparse}}$ and $z^p_{\mathrm{sparse}}$ contain the elements corresponding to sparse columns of $A^p$,
while $x^p_{\mathrm{dense}}$ and $z^p_{\mathrm{dense}}$ contain the elements corresponding to dense columns of $A^p$.
Then it defines:
\[
   M^p_{\mathrm{sparse}}
   = A^p_{\mathrm{sparse}}
     \,\operatorname{diag}(x^p_{\mathrm{sparse}})
     \,\operatorname{diag}(z^p_{\mathrm{sparse}})^{-1}
     \,(A^p_{\mathrm{sparse}})^T,
\]
\[
   U^p
   = A^p_{\mathrm{dense}}
     \,\operatorname{diag}(x^p_{\mathrm{dense}})^{\tfrac12}
     \,\operatorname{diag}(z^p_{\mathrm{dense}})^{-\tfrac12},
\]
\[ D^p = I. \]
This decomposition satisfies $M^p = M^p_{\mathrm{sparse}} + U^p D^p (U^p)^T$,
where $M^p_{\mathrm{sparse}}$ is positive definite.
In this simple case, $-(D^p)^{-1} = -I$.

\section{Other computation techniques}\label{sec:other_computation}
\subsection{Perturbation of \boldmath $M_{\text{sparse}}$ }
In practice, $M_{\text{sparse}}$ may become numerically ill-conditioned due to several factors: 
roundoff errors, rank-deficient constraint matrices, and the inherent nature of interior-point methods where approaching optimality ($\mu \to 0$) leads to increasingly poor conditioning as complementarity gaps narrow.
To ensure numerical stability of the Cholesky factorization, SDPT3 employs an adaptive perturbation strategy:
\[
  M_{\text{sparse}} \;\leftarrow\; M_{\text{sparse}} + \rho D_{M_{\text{sparse}}} + \lambda \sum_{p\in P\setminus P^{\text{u}}}A^p_{\text{sparse}} (A^p_{\text{sparse}})^T
\]
where $D_{M_{\text{sparse}}}$ is a diagonal matrix formed from the diagonal elements of $M_{\text{sparse}}$, and parameters $\rho, \lambda$ decrease geometrically with iterations.\footnote{
  Intuitively, the diagonal term $\rho D_{M_{\text{sparse}}}$ improves conditioning while preserving solution accuracy, and $\lambda \sum_{p\in P\setminus P^{\text{u}}}A^p_{\text{sparse}} (A^p_{\text{sparse}})^T$ ensures numerical stability along the constraint directions.
}
When extreme ill-conditioning is detected (e.g., condition number is greater than $10^{14}$), diagonal elements below $10^{-8}$ may be increased to $1$, but only when very few such elements exist.

These perturbations involve a delicate trade-off: while larger perturbations ensure numerical stability, they also distance the solution from the original problem, potentially compromising accuracy. SDPT3 implements a sophisticated perturbation strategy that carefully balances these competing concerns through iteration-dependent parameters and conditional application; for the complete algorithm, see the \texttt{linsysolve.m} function.

\subsection{Handling unrestricted variables}
If $\nu^p = 0$ for all $p \in P$, the problem (P) and (D) can be transformed into the following 3-parameter Homogeneous self-dual (HSD) model \cite{Wright1997}.
Given an initial point 
\[
(x_0, y_0, z_0, \tau_0, \kappa_0, \theta_0) \in 
  \operatorname{int}(\mathbb{K}) \times \mathbb{R}^m \times \operatorname{int}(\mathbb{K})
  \times \mathbb{R}^1_+ \times \mathbb{R}^1_+ \times \mathbb{R}^1_+,
\]
the HSD model is formulated as:
\[
  \begin{array}{cl}
   \min_{x,y,z,\tau,\kappa,\theta} & \bar{\alpha}\,\theta \\[3pt]
   \text{s.t.}
   & \begin{pmatrix}
       0 & -\mathcal{A} & b & -\bar{b}\\
       \mathcal{A}^T & 0 & -c & \bar{c}\\
       -b^T & c^T & 0 & -\bar{g}\\
       \bar{b}^T & -\bar{c}^T & \bar{g} & 0
     \end{pmatrix}
     \begin{pmatrix} y \\ x \\ \tau \\ \theta \end{pmatrix}
   \;+\;
     \begin{pmatrix} 0 \\ z \\ \kappa \\ 0 \end{pmatrix}
   =
     \begin{pmatrix} 0 \\ 0 \\ 0 \\ \bar{\alpha} \end{pmatrix}
  \end{array}
\]
where the parameters are defined as:
\begin{align*}
    \bar{b} &= \frac{1}{\theta_0}(b\tau_0 - \mathcal{A}x_0) \\
    \bar{c} &= \frac{1}{\theta_0}(c\tau_0 - \mathcal{A}^T y_0 - z_0) \\
    \bar{g} &= \frac{1}{\theta_0}(\inprod{c}{x_0} - b^Ty_0 + \kappa_0) \\
    \bar{\alpha} &= \frac{1}{\theta_0} (\inprod{x_0}{z_0} + \tau_0 \kappa_0)
\end{align*}

Toh et al.~\cite{toh1999} reported that if the feasible region of (P) and (D) is non-empty but lacks an interior point, converting to the HSD model yields superior accuracy despite increased computational cost.
This is particularly important when unrestricted blocks exist (i.e., $\exists p\in P$ such that $\mathbb{K}^p=\mathbb{R}^{n_p}$), as the dual cone $(\mathbb{K}^p)^*=\{0\}$ implies that there is no interior point in the dual feasible region.
SDPT3 applies this HSD transformation automatically when $\nu^p = 0$ for all $p \in P$ and $\mathbb{K}^p=\mathbb{R}^{n_p}$ for some $p\in P$.
The standard algorithm cannot be directly applied to the HSD model and requires modifications, which are detailed in \cite{toh1999}.

\medskip

When $\nu^p \neq 0$ for some $p \in P$, the problem (P) and (D) cannot be converted to the HSD model.
As seen from \eqref{eq:Schur_complement_Mat}, unrestricted blocks (i.e., $\mathbb{K}^p = \mathbb{R}^{n_p}$) require solving augumented systems involving $A^{\text{u}}$, which increases computational cost.
To avoid this overhead, SDPT3 eliminates unrestricted variables $x^p \in \mathbb{R}^{n_p}$ by splitting them into non-negative components:
\[
   x^p = x^p_+ - x^p_-, \quad \text{where } (x^p_+, x^p_-) \in \mathbb{R}^{2n_p}_+
\]

This transformation, while computationally efficient, introduces numerical challenges.
As iterations progress, both $x^p_+$ and $x^p_-$ typically grow large while their dual counterparts $z^p_+$ and $z^p_-$ become small, causing the complementarity products $\operatorname{diag}(x^p_\pm) \operatorname{diag}(z^p_\pm)$ to become severely ill-conditioned.

To mitigate these issues, SDPT3 employs heuristic stabilization techniques.
At each iteration, the primal variables are recentered:
\[
   x^p_+ \leftarrow x^p_+ - 0.8 \min(x^p_+, x^p_-), \quad
   x^p_- \leftarrow x^p_- - 0.8 \min(x^p_+, x^p_-)
\]
while the dual variables are perturbed by adding $\alpha \mu e^p$, where $e^p$ is the identity element defined in Section 3.1 and $\alpha = 0.1$ is an adaptively defined small constant, to prevent them from approaching zero too rapidly.

\medskip

Similar numerical issues arise even in linear blocks ($\mathbb{K}^p=\mathbb{R}^{n_p}_+$) when pairs of variables effectively represent a single unrestricted variable.
Specifically, if variables $x^p_i$ and $x^p_j$ satisfy $(a^p_k)_i = -(a^p_k)_j$ for all constraints $k=1,\ldots,m$ and $(c^p)_i = -(c^p)_j$ with $\nu^p = 0$, they exhibit the same instability pattern.
SDPT3 detects such implicit unrestricted variables and applies the same stabilization techniques described above.

\subsection{Preprocessing for model transformation}
\subsubsection{Transforming complex semidefinite variables to real}
In applications such as control engineering, optimization problems on complex semidefinite cones may arise, but they can be reduced to problems on real semidefinite cones.

Define the following:
\begin{itemize}
    \item Complex matrix space $\mathbb{C}^{m\times n}$
    \item Set of Hermitian matrices $\mathbb{H}^n=\{a \in \mathbb{C}^{n\times n} \mid a = a^H\}$ where $a^H$ denotes the conjugate transpose
    \item Complex positive semidefinite cone $\mathbb{H}^n_+=\{a \in \mathbb{H}^n \mid z^H a z \geq 0 ~(\forall z\in \mathbb{C}^n)\}$
    \item $\bar{\mathbb{S}}^{2n}_+ = \left\{\left(\begin{smallmatrix}
      R & -S\\
      S & R
  \end{smallmatrix} \right) \in \mathbb{S}^{2n}_+ \;\middle|\; R\in \mathbb{S}^n, ~ S^T=-S\right\}.$
\end{itemize}
Any Hermitian matrix $a \in \mathbb{H}^n$ can be uniquely decomposed as $a = R + iS$ where $R = \operatorname{real}(a) \in \mathbb{S}^n$ is the real part which is a real symmetric matrix and $S = \operatorname{imag}(a) \in \mathbb{R}^{n\times n}$ is the imaginary part which is a real skew-symmetric matrix, i.e., $S^T = -S$.

We define the real embedding $\Gamma: \mathbb{H}^n\to \mathbb{S}^{2n}$ as
\[
  \Gamma(a) 
  = \begin{pmatrix}
       \operatorname{real}(a) & -\operatorname{imag}(a) \\
       \operatorname{imag}(a) & \operatorname{real}(a)
     \end{pmatrix}.
\]
Then the mapping $\Gamma$ establishes a cone isomorphism between $\mathbb{H}^n_+$ and $\bar{\mathbb{S}}^{2n}_+$:
\[
  a\in \mathbb{H}^n_+
   \;\;\Longleftrightarrow\;\;
  \Gamma(a)\in \bar{\mathbb{S}}^{2n}_+.
\]
This equivalence follows from the identity: for any $u,v\in\mathbb{R}^n$ and $z=u+iv\in\mathbb{C}^n$,
\[
  z^H a z = 
  \begin{pmatrix}u\\ v\end{pmatrix}^T \Gamma(a) \begin{pmatrix}u\\ v\end{pmatrix}.
\]
The left-hand side is non-negative for all $z \in \mathbb{C}^n$ if and only if the right-hand side is non-negative for all $(\begin{smallmatrix}u\\v\end{smallmatrix}) \in \mathbb{R}^{2n}$.

Thus, by transforming complex variables $x \in \mathbb{H}^{n_p}_+$ to real variables $\bar{x} = \Gamma(x) \in \bar{\mathbb{S}}^{2n_p}_+$, problems on complex semidefinite cones can be reduced to problems on real semidefinite cones.

\begin{example}
\begin{equation*}
    \left|
    \begin{array}{cl}
        \min & \inprod{\begin{pmatrix}
            2 & 1-i \\
            1+i & 3
        \end{pmatrix}}{x} \\
        s.t. 
        & \inprod{\begin{pmatrix}
            1 & 0 \\ 
            0 & 1
        \end{pmatrix}}{x} = 4 \\
        & x\in \mathbb{H}^2_+ 
    \end{array}
    \right.
    \Longleftrightarrow
    \left|
    \begin{array}{cl}
        \min & \inprod{\begin{pmatrix}
            2 & 1 & 0 & 1 \\
            1 & 3 & -1 & 0 \\
            0 & -1 & 2 & 1 \\
            1 & 0 & 1 & 3
        \end{pmatrix}}{\bar{x}} \\
        s.t. 
        & \inprod{\begin{pmatrix}
            1 & 0 & 0 & 0 \\ 
            0 & 1 & 0 & 0 \\ 
            0 & 0 & 1 & 0 \\ 
            0 & 0 & 0 & 1
        \end{pmatrix}}{\bar{x}} = 4 \\
        & \bar{x}\in \bar{\mathbb{S}}^4_+ = \left\{\begin{pmatrix}
            R & -S\\
            S & R
        \end{pmatrix} \in \mathbb{S}^4_+ \;\middle|\; R\in \mathbb{S}^2, ~ S^T=-S\right\}
    \end{array}
    \right.
\end{equation*}
\end{example}

\subsubsection{Converting diagonal blocks to linear variables}
The condition $x\in \mathbb{S}^1_+$ is equivalent to $x\in \mathbb{R}^1_+$, since a $1\times 1$ positive semidefinite matrix is simply a non-negative scalar.
In this case, treating $x$ as a non-negative real variable improves the computational efficiency of the interior-point method.
Furthermore, if the variable matrix in $\mathbb{S}^{n_p}_+$ contains isolated diagonal elements (i.e., elements that do not interact with off-diagonal entries in the constraints), converting them to non-negative real variables yields similar computational efficiency improvements.

Specifically, for an integer $i$, if the $i$-th diagonal element is isolated, i.e.,
$(c^p)_{ij}=(c^p)_{ji}=0$ and $(a^p_k)_{ij}=(a^p_k)_{ji}=0$ for all $j\neq i$ and all $k$,
SDPT3 transforms $x\in \mathbb{S}^{n_p}_+$ to $(\bar{x},\,\hat{x})\in \mathbb{S}^{n_p-1}_+\times \mathbb{R}^1_+$,
where $\bar{x}$ is the $(n_p-1)\times(n_p-1)$ matrix obtained by removing the $i$-th row and column from $x$,
and $\hat{x} = x_{ii} \geq 0$ represents the isolated diagonal element.

\begin{example}
\begin{equation*}
    \left|
    \begin{array}{cl}
        \min & \inprod{\begin{pmatrix}
            3 & 0 & 1 \\
            0 & 5 & 0 \\
            1 & 0 & 2
        \end{pmatrix}}{x} \\
        s.t. 
        & \inprod{\begin{pmatrix}
            1 & 0 & 0 \\ 
            0 & 2 & 0 \\ 
            0 & 0 & 3
        \end{pmatrix}}{x} = 1 \\
        & x\in \mathbb{S}^3_+ 
    \end{array}
    \right.
    \Longleftrightarrow
    \left|
    \begin{array}{cl}
        \min & \inprod{\begin{pmatrix}
            3 & 1 \\
            1 & 2
        \end{pmatrix}}{\bar{x}} + 5\hat{x} \\
        s.t. 
        & \inprod{\begin{pmatrix}
            1 & 0 \\ 0 & 3
        \end{pmatrix}}{\bar{x}} + 2\hat{x} = 1 \\
        & \bar{x}\in \mathbb{S}^2_+, \quad \hat{x} \in \mathbb{R}^1_+
    \end{array}
    \right.
\end{equation*}
\end{example}
Indeed, for example, the instance \texttt{qpG11} from SDPLIB contains many such diagonal blocks that benefit from this transformation technique.

\subsubsection{Ensuring algorithmic assumptions through variable augmentation}
The interior-point method algorithm introduced in this paper requires $m \geq 1$ and $\exists p\in P ~ \text{s.t.} ~ \nu^p=0$ (mainly for \eqref{eq:NewtonKKT} and \eqref{eq:mu}).
When these prerequisites are not satisfied, SDPT3 adds one artificial non-negative variable $x^{p_{\max} + 1}$ and the corresponding constraint
\[
  -\sum_{p\in P} \inprod{e^p}{x^p}_p + x^{p_{\max} + 1} = 0
\]
to ensure the algorithmic assumptions hold.
Specifically, this augmentation introduces the following additional problem parameters:
\begin{itemize}
    \item $\mathbb{K}^{p_{\max} + 1} = \mathbb{R}^1_+$
    \item $a^{p}_{m+1}=e^p \;\; (p \in P) \;\;$ and $\;\; a^{p_{\max} + 1}_{m+1} = 1$
    \item $b_{m+1} = 0$
    \item $c^{p_{\max} +1}=0$
    \item $\nu^{p_{\max} + 1} = 0$
\end{itemize}
Note that when $m = 0$ (no existing constraints), this creates the first constraint with $b_1 = 0$. When $m \geq 1$, this extends the existing constraint vector $b$ by appending a new element $b_{m+1} = 0$.

\subsubsection{Reordering matrices for Cholesky efficiency}
To verify the positive definiteness of variables $x^p, z^p$, SDPT3 uses Cholesky decomposition, which provides a numerically stable test.
Since the interior-point method performs Cholesky decomposition at each iteration, minimizing fill-in is crucial for computational efficiency.
SDPT3 reorders matrix variables to minimize fill-in during Cholesky decomposition.
To identify the sparsity pattern, SDPT3 first constructs the aggregate matrix
\[
  t^p = |c^p| + \sum_{k=1}^m |a^p_k|
\]
which represents the combined sparsity structure of all coefficient matrices.
The Reverse Cuthill-McKee algorithm \cite{Cuthill1969,Golub2013} is then applied to obtain the permutation $\sigma$ that reduces the bandwidth (the maximum distance of non-zero elements from the diagonal).
Since reducing bandwidth typically reduces fill-in during Cholesky decomposition, this reordering improves computational efficiency.\footnote{
  While the benefits are problem-dependent, this reordering often improves efficiency in practice. 
  The Newton system matrix $M = \mathcal{A}\mathcal{H}\mathcal{A}^T$ tends to inherit sparsity from the coefficient matrices $\mathcal{A}$, and iterates $x^p, z^p$ often reflect the coefficient structure $(a_k^p)$, particularly near optimality. 
  The reduced bandwidth typically leads to less fill-in during Cholesky decompositions throughout the algorithm.}
The variables are reordered according to $(\bar{x}^p)_{ij} = (x^p)_{\sigma(i)\sigma(j)}$.

\begin{example}
\[
    t^p = \begin{pmatrix}
        3 & 0 & 1 \\\
        0 & 5 & 0 \\\
        1 & 0 & 2
    \end{pmatrix}
\]
Applying the Reverse Cuthill-McKee algorithm yields $\sigma(1)=3, ~ \sigma(2)=1, ~ \sigma(3)=2$.
Transforming variables using this permutation results in $\bar{t}^p$:
\[
    \bar{t}^p = \begin{pmatrix}
        3 & 1 & 0 \\
        1 & 2 & 0 \\
        0 & 0 & 5
    \end{pmatrix}
\]
which reduces the bandwidth from 2 to 1.
\end{example}

\subsubsection{Utilizing matrix symmetry}
In blocks where $\mathbb{K}^p=\mathbb{S}_+^p$, exploiting the symmetry of matrices reduces memory usage and computation time for inner products.
Given a real symmetric matrix $a\in \mathbb{S}^n$, the symmetric vectorization operator $\operatorname{svec}: \mathbb{S}^n \rightarrow \mathbb{R}^{n(n+1)/2}$ extracts the upper triangular part column by column, scaling off-diagonal elements by $\sqrt{2}$ to preserve inner products:
\[ \operatorname{svec}(a) = (a_{11}, \quad \sqrt{2}\,a_{12},\, a_{22},\, \quad \sqrt{2}\,a_{13},\, \sqrt{2}\,a_{23},\, a_{33},\, \quad \ldots \quad  \sqrt{2}\,a_{1n},\, \sqrt{2}\,a_{2n},\, \ldots,\, a_{nn})^T \]
\begin{example}
\[a = \begin{pmatrix}
    1 & 2 & 3\\
    2 & 4 & 5\\
    3 & 5 & 6
\end{pmatrix} \rightarrow \operatorname{svec}(a) = (1, ~ 2\sqrt{2}, ~ 4, ~ 3\sqrt{2}, ~ 5\sqrt{2}, ~ 6)^T\]
\end{example}

\medskip

Indeed, the following equality holds for any symmetric matrices $a, b \in \mathbb{S}^n$:
\[\inprod{a}{b} = \operatorname{svec}(a)^T \operatorname{svec}(b).\]
Therefore, storing $\operatorname{svec}(a)$ instead of $a$ halves memory usage while making inner product computations more efficient.
SDPT3 stores coefficient matrices $a^p_k$ in the $\operatorname{svec}$ format for efficiency, while maintaining $c^p$, $x^p$, and $z^p$ in matrix form for algorithmic operations.

\section{Summary}
The pseudo-code for the primal-dual path-following interior-point method incorporating various techniques discussed in the preceding sections is presented as Algorithm~\ref{alg:pd_ipm}. 
\begin{algorithm}
\caption{Primal-dual path-following interior-point method}
\label{alg:pd_ipm}
\begin{algorithmic}[1]
\STATE \textbf{Preprocessing:} Transform problem as described in Section~\ref{sec:other_computation}. 
Generate infeasible starting point $(x,y,z) \leftarrow (x^0, y^0, z^0)$
\FOR{$k = 1, 2, \ldots $}
    \STATE Compute primal and dual residuals $R_{\text{prim}}, R_{\text{dual}}^p$, and check convergence
    \STATE Compute current duality gap $\mu$ and assess stopping criteria
    \IF{convergence achieved}
        \RETURN $(x, y, z)$
    \ENDIF
    
    \STATE \textbf{System setup:}
    \STATE Compute coefficient matrices $M_{\text{sparse}}$, $A^{\text{u}}$, $U$, $-D^{-1}$ for augmented system
    \STATE Apply matrix perturbation $M_{\text{sparse}} \leftarrow M_{\text{sparse}} + \rho D_{M_{\text{sparse}}} + \lambda \sum A^p_{\text{sparse}}(A^p_{\text{sparse}})^T$
    \STATE Compute Cholesky decomposition of $M_{\text{sparse}}$ (or LU if Cholesky fails)
    \STATE Compute $\mathcal{A}\mathcal{H}R_{\text{dual}}$ and $R^{\text{u}}_{\text{dual}}$
    \STATE 
    \STATE \textbf{Predictor step:}
    \STATE Compute complementarity residual $R_{\text{comp}}^{p,\text{pred}}$
    \STATE Compute $\mathcal{A}\mathcal{E}^{-1}R_{\text{comp}}^{\text{pred}}$ 
    \STATE Compute right-hand side $h_{\text{pred}} = R_{\text{prim}} + \mathcal{A}\mathcal{E}^{-1}R_{\text{comp}}^{\text{pred}} - \mathcal{A}\mathcal{H}R_{\text{dual}}$
    \STATE Solve augmented system \eqref{eq:Schur_complement_Mat_aug} to obtain $\delta y, \delta x^{\text{u}}$
    \STATE Compute $\delta x, \delta z$ from equations \eqref{eq:sol_x_z}
    \STATE Compute step lengths $\alpha_{P}, \alpha_{D}$
    \STATE Compute centering parameter $\sigma$ with adaptive exponent $\psi$
    \STATE
    \STATE \textbf{Corrector step:}
    \STATE Compute complementarity residual $R_{\text{comp}}^{p,\text{corr}}$
    \STATE Compute $\mathcal{A}\mathcal{E}^{-1}R_{\text{comp}}^{\text{corr}}$
    \STATE Compute right-hand side $h_{\text{corr}} = R_{\text{prim}} + \mathcal{A}\mathcal{E}^{-1}R_{\text{comp}}^{\text{corr}} - \mathcal{A}\mathcal{H}R_{\text{dual}}$
    \STATE Solve augmented system \eqref{eq:Schur_complement_Mat_aug} to obtain $\Delta y, \Delta x^{\text{u}}$
    \STATE Compute $\Delta x, \Delta z$ from equations \eqref{eq:sol_x_z}
    \STATE Compute step lengths $\beta_P, \beta_D$ ensuring positive definiteness
    \STATE Update: $x \leftarrow x + \beta_P \Delta x$, $y \leftarrow y + \beta_D \Delta y$, $z \leftarrow z + \beta_D \Delta z$

    \STATE 
    \IF{heuristic stopping criteria met}
        \RETURN $(x, y, z)$
    \ENDIF
    \STATE Apply stabilization for free variables
\ENDFOR
\end{algorithmic}
\end{algorithm}

The algorithm's practical success stems from the careful integration of techniques presented throughout this paper: the predictor-corrector framework that balances aggressive progress with numerical stability; sophisticated linear system solving with adaptive factorization and iterative refinement; carefully designed initial points; systematic sparsity exploitation; adaptive perturbation strategies for handling ill-conditioning, etc.

Furthermore, the algorithm achieves computational efficiency by exploiting the fact that the coefficient matrices $M_{\text{sparse}}$, $A^{\text{u}}$, $U$, and $-D^{-1}$ remain unchanged between predictor and corrector steps. This invariance allows the reuse of the expensive Cholesky (or LU) factorization of the augmented system across both linear system solves. The only modification required between steps involves the right-hand side vectors, specifically the complementarity residuals: $R_{\text{comp}}^{p,\text{pred}}$ for the predictor versus $R_{\text{comp}}^{p,\text{corr}}$ for the corrector, where the latter incorporates both centering terms ($\sigma \mu$) and second-order corrections derived from the predictor step.

Note that the actual SDPT3 implementation is highly optimized with numerous computational enhancements that make the practical code significantly more complex than this simplified presentation. 
Especially, the efficient handling of low-rank structured data, the 3-parameter homogeneous self-dual model for problems with free variables, and sophisticated heuristic stopping criteria are not fully presented in this paper. 
Readers interested in the complete implementation details are encouraged to consult the original SDPT3 papers \cite{toh1999,Toh2012,tutuncu2003} and the source code.

\section*{Acknowledgments}
The author would like to thank Professor Kim-Chuan Toh for his valuable responses to inquiries
regarding the SDPT3 algorithm and implementation details.

\clearpage
\phantomsection
\addcontentsline{toc}{part}{Appendix}
\begin{center}
{\Large\bfseries\textsf{Appendix}}
\end{center}
\vspace{2em}

\appendix
\section{A Guide to the Derivation of Equations \eqref{eq:Einv_Rcomp}--\eqref{eq:NT_M}} \label{sec:guide_for_dir_eq}

\paragraph{Semidefinite Cone Case: $\mathbb{K}^p=\mathbb{S}_+^{n_p}$}
In this case, the derivation for the HKM direction is straightforward. Note that $G^p=(z^p)^\frac{1}{2}$ and 
\begin{equation*}
\mathcal{E}\Delta x = ((z^p)^\frac{1}{2}\Delta x^p)\circ ((z^p)^{-\frac{1}{2}}z^p)=(z^p)^\frac{1}{2} \Delta x^p (z^p)^\frac{1}{2}.
\end{equation*}
For the NT direction, the derivation is non-trivial and can be found in \cite{todd1998}.

\paragraph{Second-Order Cone Case: $\mathbb{K}^p=\mathbb{Q}^{n_p}$}
We define the linear operator $\operatorname{Arw}: \mathbb{R}^{n_p} \rightarrow \mathbb{R}^{n\times n}$ as follows. For $f\in \mathbb{R}^{n^p}$,
\begin{equation*}
    \operatorname{Arw}(f) = 
    \left(
    \begin{array}{cc}
        f_0 & \bar{f}^T \\
        \bar{f} & f_0 I
    \end{array}
    \right)
\end{equation*}
Note that $x^p \circ z^p = \operatorname{Arw}(x^p)\,z^p$. Therefore, the following expressions hold:
\begin{align*}
    (\mathcal{E}^p)^{-1}R^p_{comp} &= (G^p)^{-1} \cdot \operatorname{Arw}\big((G^p)^{-1} z^p\big)^{-1} \cdot R^p_{comp}\\
    \mathcal{H}^p R^p_{dual} &= (G^p)^{-1} \cdot \operatorname{Arw}\big((G^p)^{-1} z^p\big)^{-1} \cdot \operatorname{Arw}\big(G^p x^p\big) \cdot(G^p)^{-1} \cdot R^p_{dual}\\
    \mathcal{A}^p\mathcal{H}^p(\mathcal{A}^p)^T &= A^p \cdot (G^p)^{-1} \cdot \operatorname{Arw}\big((G^p)^{-1} z^p\big)^{-1} \cdot \operatorname{Arw}\big(G^p x^p\big) \cdot(G^p)^{-1} \cdot (A^p)^T
\end{align*}
These relations provide the general form for the required computations. 
Furthermore,
\begin{equation*}
    \operatorname{Arw}(f)^{-1} = \frac{1}{\gamma(f)^2} \left(\begin{array}{cc}
        f_0 & -\bar{f}^T \\
        -\bar{f} & \frac{1}{f_0}(\gamma(f)^2 I + \bar{f}\bar{f}^T)
    \end{array}\right),
    \qquad
    (G^p)^{-1} = \frac{1}{\omega^p} \begin{pmatrix}
        t^p_0 & -(\bar{t}^p)^T \\
        -\bar{t}^p & I+\frac{1}{1 + t^p_0} \bar{t}^p(\bar{t}^p)^T
    \end{pmatrix}
\end{equation*}
where we assume $\gamma(t^p)=1$. 
Equations \eqref{eq:HKM_HRd}--\eqref{eq:NT_M} can be derived by substituting the specific form of $G^p$ and expanding these expressions algebraically.

\paragraph{Linear Cone Case: $\mathbb{K}^p=\mathbb{R}^{n_p}_+$}
In this case, we always use $G^p=I$, hence the derivation is straightforward by noting that $\mathcal{E}^p \Delta x^p = \operatorname{diag}(z^p) \Delta x^p$.

\section{SDPT3 Input Data Format}
This section provides a detailed explanation of the input data format for SDPT3.
Since SDPT3 is implemented in MATLAB, we use MATLAB notation throughout this section:
\begin{itemize}
\item \texttt{[a,b]}: Horizontal concatenation of vectors or matrices \texttt{a} and \texttt{b}, i.e., $\begin{pmatrix} a \; b \end{pmatrix}$.
\item \texttt{[a;b]}: Vertical concatenation of \texttt{a} and \texttt{b}, i.e., $\begin{pmatrix} a \\ b \end{pmatrix}$.
\end{itemize}

SDPT3 is called as follows:
\begin{lstlisting}[language=Matlab]
[x, y, info] = sdpt3(blk, At, C, b, OPTIONS);
\end{lstlisting}

For problems (P) and (D), the SDPT3 input arguments are set as follows.

\begin{itemize}
    \item \texttt{blk} is a $p_{\max} \times 2$ cell array\footnote{MATLAB's cell array is a container that can hold different data types, similar to Python's lists.} that specifies the type and dimension of each subcone $\mathbb{K}^p$: 
\begin{lstlisting}[language=Matlab,escapechar=\@]
blk{p, 1} = @$\texttt{type}^p$@;  blk{p, 2} = @$n_p$@;
\end{lstlisting}
where the type identifier for each cone is defined as:
\[
\texttt{type}^p = \begin{cases}
\texttt{'s'} & \text{if } \mathbb{K}^p = \mathbb{S}_+^{n_p}, \\
\texttt{'q'} & \text{if } \mathbb{K}^p = \mathbb{Q}^{n_p}, \\
\texttt{'l'} & \text{if } \mathbb{K}^p = \mathbb{R}^{n_p}_+, \\
\texttt{'u'} & \text{if } \mathbb{K}^p = \mathbb{R}^{n_p}. 
\end{cases}
\]
\item \texttt{At} is a $p_{\max} \times 1$ cell array which specifies the constraint matrices for each subcone: For semidefinite cones, 
\begin{lstlisting}[language=Matlab,escapechar=\@]
At{p} = [@$\operatorname{svec}(a^p_1), \operatorname{svec}(a^p_2), ..., \operatorname{svec}(a^p_m)$@];
\end{lstlisting}
For other cones,
\begin{lstlisting}[language=Matlab,escapechar=\@]
At{p} = [@$a^p_1, a^p_2, ..., a^p_m$@];
\end{lstlisting}
\item \texttt{C} is a $p_{\max} \times 1$ cell array that stores the objective function coefficients for each subcone:
\begin{lstlisting}[language=Matlab,escapechar=\@]
  C{p} = @$c^p$@;
\end{lstlisting}
\item \texttt{b} is an $m \times 1$ vector representing the right-hand side of equality constraints:
\begin{lstlisting}[language=Matlab,escapechar=\@]
  b = [@$b_1; b_2; ...; b_m$@];
\end{lstlisting}
\item \texttt{OPTIONS} is a structure that allows various option settings. 
The coefficient $\nu^p$ can be set through \texttt{OPTIONS.parbarrier} which is an optional $p_{\max} \times 1$ vector:
\begin{lstlisting}[language=Matlab,escapechar=\@]
OPTIONS.parbarrier = [@$\nu^1; \nu^2; ...; \nu^{p_{\max}}$@];
\end{lstlisting}
\end{itemize}

\begin{example}
The following optimization problem demonstrates how to provide input to SDPT3:
\[
\left|
\begin{array}{cl}
\max  & 6y_1 + 4y_2 + 5y_3 \\
\mathrm{s.t.} 
& 16y_1 - 14y_2 + 5y_3 \leq -3 \\
& 7y_1 + 2y_2 \leq 5 \\
& \left\| \begin{pmatrix} 8y_1 + 13y_2 - 12y_3 - 2 \\ -8y_1 + 18y_2 + 6y_3 - 14 \\ y_1 - 3y_2 - 17y_3 - 13 \end{pmatrix} \right\|  \leq -24y_1 - 7y_2 + 15y_3 + 12 \\
& \left\| \begin{pmatrix} y_1 \\ y_2 \\ y_3 \end{pmatrix} \right\|  \leq 10 \\
& \begin{pmatrix} 
7y_1 + 3y_2 + 9y_3 & -5y_1 + 13y_2 + 6y_3 & y_1 - 6y_2 - 6y_3 \\
-5y_1 + 13y_2 + 6y_3 & y_1 + 12y_2 - 7y_3 & -7y_1 - 10y_2 - 7y_3 \\
y_1 - 6y_2 - 6y_3 & -7y_1 - 10y_2 - 7y_3 & -4y_1 - 28y_2 - 11y_3
\end{pmatrix}  \preceq \begin{pmatrix} 68 & -30 & -19 \\ -30 & 99 & 23 \\ -19 & 23 & 10 \end{pmatrix}
\end{array}
\right.
\]
This problem can be categorized as a dual problem (D) with:
\begin{itemize}
    \item $m = 3$, 
    $\quad b = \begin{pmatrix} 6,  4,  5 \end{pmatrix}^T$
    \item $p_{\max} = 4$, 
    $\quad \nu^p = 0$ for all $p=1,\ldots,4$
    \item $\mathbb{K}^1 = \mathbb{R}^{n_1}_+$, 
    $\mathbb{K}^2 = \mathbb{Q}^{n_2}$, 
    $\mathbb{K}^3 = \mathbb{Q}^{n_3}$, 
    $\mathbb{K}^4 = \mathbb{S}^{n_4}_+$
    \item $n_1 = 2$, 
    $n_2 = 4$, 
    $n_3 = 4$, 
    $n_4 = 3$
    \item $c^1 = \begin{pmatrix} -3 \\ 5 \end{pmatrix}$, 
    $a^1_1 = \begin{pmatrix} 16 \\ 7 \end{pmatrix}$, 
    $a^1_2 = \begin{pmatrix} -14 \\ 2 \end{pmatrix}$, 
    $a^1_3 = \begin{pmatrix} 5 \\ 0 \end{pmatrix}$
    \item $c^2 = \begin{pmatrix} 12 \\ -2 \\ -14 \\ -13 \end{pmatrix}$, 
    $a^2_1 = \begin{pmatrix} 24 \\ -8 \\ 8 \\ -1 \end{pmatrix}$, 
    $a^2_2 = \begin{pmatrix} 7 \\ -13 \\ -18 \\ 3 \end{pmatrix}$, 
    $a^2_3 = \begin{pmatrix} -15 \\ 12 \\ -6 \\ 17 \end{pmatrix}$ 
    \item $c^3 = \begin{pmatrix} 10 \\ 0 \\ 0 \\ 0 \end{pmatrix}$, 
    $a^3_1 = \begin{pmatrix} 0 \\ -1 \\ 0 \\ 0 \end{pmatrix}$, 
    $a^3_2 = \begin{pmatrix} 0 \\ 0 \\ -1 \\ 0 \end{pmatrix}$, 
    $a^3_3 = \begin{pmatrix} 0 \\ 0 \\ 0 \\ -1 \end{pmatrix}$
    \item $c^4 = \begin{pmatrix} 68 & -30 & -19 \\ -30 & 99 & 23 \\ -19 & 23 & 10 \end{pmatrix}$, \\
    $a^4_1 = \begin{pmatrix} 7 & -5 & 1 \\ -5 & 1 & -7 \\ 1 & -7 & -4 \end{pmatrix}$, 
    $a^4_2 = \begin{pmatrix} 3 & 13 & -6 \\ 13 & 12 & -10 \\ -6 & -10 & -28 \end{pmatrix}$, 
    $a^4_3 = \begin{pmatrix} 9 & 6 & -6 \\ 6 & -7 & -7 \\ -6 & -7 & -11 \end{pmatrix}$
\end{itemize}

The SDPT3 input becomes:
\begin{lstlisting}[language=Matlab]
blk{1,1} = 'l'; blk{1,2} = 2;
blk{2,1} = 'q'; blk{2,2} = 4;
blk{3,1} = 'q'; blk{3,2} = 4;
blk{4,1} = 's'; blk{4,2} = 3;

At{1} = [16, -14, 5;
         7,   2, 0];

At{2} = [24,   7, -15;
         -8, -13,  12;
          8, -18,  -6;
         -1,   3,  17];

At{3} = [ 0,  0,  0;
         -1,  0,  0;
          0, -1,  0;
          0,  0, -1];

A1_sdp = [7, -5, 1; -5, 1, -7; 1, -7, -4];
A2_sdp = [3, 13, -6; 13, 12, -10; -6, -10, -28];
A3_sdp = [9, 6, -6; 6, -7, -7; -6, -7, -11];
pblk = {'s', 3};
At{4} = [svec(pblk,A1_sdp), svec(pblk,A2_sdp), svec(pblk,A3_sdp)];
%% that is equivalent to:
% At(4) = svec(pblk, {A1_sdp, A2_sdp, A3_sdp});
%% or explicitly:
% s = sqrt(2);
% At{4} = [ 7  ,   3  ,   9  ;
%          -5*s,  13*s,   6*s;
%           1  ,  12  ,  -7  ;
%           1*s,  -6*s,  -6*s;
%          -7*s, -10*s,  -7*s;
%          -4  , -28  , -11  ];

C{1} = [-3; 5];
C{2} = [12; -2; -14; -13];
C{3} = [10; 0; 0; 0];
C{4} = [68, -30, -19; -30, 99, 23; -19, 23, 10];

b = [6; 4; 5];

[x, y, info] = sdpt3(blk, At, C, b);
\end{lstlisting}
\end{example}

\paragraph{Efficient Usage Note}
In SDPT3, it is possible to group multiple semidefinite cones or second-order cones into a single cell entry, rather than treating each cone as a separate block. 
Since MATLAB's iteration over cell arrays can be slow, such consolidated input can be computationally more efficient, particularly when dealing with numerous small cones of the same type. 
For detailed information on this advanced usage, see \cite{Toh2012}.

\bibliographystyle{unsrtnat}
\bibliography{references}
\end{document}